\documentclass[11pt,a4paper]{article}
\usepackage[T1]{fontenc}
\usepackage[utf8]{inputenc}
\usepackage{lmodern}
\usepackage[margin=28mm,headheight=14pt]{geometry}
\usepackage{amsmath,amssymb,amsthm,mathtools,mathrsfs}
\usepackage{enumitem,graphicx,xcolor,booktabs,microtype}
\usepackage{hyperref}
\usepackage{bookmark}
\hypersetup{colorlinks=true,linkcolor=blue!45!black,citecolor=blue!45!black,urlcolor=blue!45!black,
 pdftitle={Leray’s Problem for the Stationary Navier–Stokes Equations with Low-Regularity Data in Planar and Three-Dimensional Axisymmetric Lipschitz Domains},
 pdfsubject={Stationary Navier--Stokes equations; solenoidal energy dual; pressure representatives},
 pdfkeywords={Leray problem, Lipschitz domain, energy trace, solenoidal dual, pressure}}
\numberwithin{equation}{section}
\newtheorem{theorem}{Theorem}[section]
\newtheorem{lemma}[theorem]{Lemma}
\newtheorem{proposition}[theorem]{Proposition}

\theoremstyle{definition}

\theoremstyle{remark}
\newtheorem{remark}[theorem]{Remark}

\setlist{topsep=4pt,itemsep=2pt,parsep=0pt}
\allowdisplaybreaks[2]
\title{Leray’s Problem for the Stationary Navier-Stokes Equations with Low-Regularity Data in Planar and Three-Dimensional Axisymmetric Lipschitz Domains}
\author{Zhibing Zhang\\[0.3em]
  {\small School of Mathematics and Statistics, Anqing Normal University}\\
  {\small  Anqing 246133, China}\\[0.3em]
  {\small \texttt{zhibingzhang29@126.com}} }
\date{}
\begin{document}
\maketitle

\begin{abstract}
In this paper, we study the steady Navier-Stokes Dirichlet problem in bounded planar
and three-dimensional axisymmetric domains with finitely many disjoint Lipschitz
boundary components. For prescribed velocities in
$H^{1/2}(\partial\Omega)$ and forces in the dual of the solenoidal energy
space $V=H_0^1(\operatorname{div}0,\Omega)$, we establish finite-energy
solvability under the sole compatibility condition of zero total
boundary flux. Individual component fluxes are subject to no smallness
or sign restrictions. In three dimensions the domain and data are
axisymmetric, and swirl is allowed. The result addresses the regularity
of the domain and data within these two geometric classes. The proof
combines Leray-Schauder continuation with pressure compactness on
Lipschitz boundary collars and a total-head level-set argument.
Continuation of the limiting Euler fields through complementary
components supplies separating cycles, while viscous identities are
evaluated inside the physical domain. A nonlinear correction of smooth
force approximations controls the volume and surface errors in the
passage to energy-dual data. An isometric extension from $V'$ to
$H^{-1}$ preserves the prescribed symmetries and fixes the force
representative used in the compactness argument. For each selected
$H^{-1}$ representative, the resulting velocity admits an associated
mean-zero $L^2$ pressure; the velocity equation depends only on the force
class modulo gradients.
\end{abstract}

\tableofcontents

\section{Introduction and main results}
\label{sec:main}

\subsection{The Leray problem and its geometric setting}

The nonhomogeneous Dirichlet problem for a steady incompressible viscous
flow asks for a velocity $u$ and a pressure $P$ satisfying
\[
 \begin{cases}
 -\nu\Delta u+(u\cdot\nabla)u+\nabla P=F & \text{in }\Omega,\\
 \operatorname{div}u=0 & \text{in }\Omega,\\
 \operatorname{Tr}u=a & \text{on }\partial\Omega,
 \end{cases}
 \qquad \nu>0.
\]
Here $F$ denotes a force acting on general test fields; below we
distinguish such a representative from its restriction to solenoidal
fields. On a multiply connected domain, incompressibility imposes a
condition on the sum of the fluxes through the boundary components.
The Leray problem is to obtain existence under this necessary
condition without imposing zero flux on each component separately.
The classical theory, its flux restrictions, and the later general-flux
results are presented in \cite{Galdi,KPRbook}.

Let $\Omega\subset\mathbb R^n$ be bounded and connected, with disjoint
connected Lipschitz boundary components $\Gamma_0,\ldots,\Gamma_N$.
We assume the finite outer-boundary and hole geometry specified in
Section~\ref{sec:continuation}; $\Gamma_0$ is the outer boundary.
The normal $n_\Omega$ always points out of the fluid, including on the
boundaries of the holes. For a prescribed trace
$a\in H^{1/2}(\partial\Omega;\mathbb R^n)$, put
\begin{equation}\label{eq:1.flux}
 F_j=\int_{\Gamma_j}a\cdot n_\Omega\,dS,
 \qquad \sum_{j=0}^N F_j=0.
\end{equation}
Only the sum in \eqref{eq:1.flux} is imposed by incompressibility.
All boundary Sobolev spaces are the ordinary Lipschitz trace spaces.

\subsection{Earlier results and the scope of the present work}
\label{subsec:literature}

Korobkov, Pileckas and Russo
\cite[Theorems~1.1 and~4.1]{KPR2015} established solvability under
the sole total-flux condition in bounded planar domains and in the
three-dimensional axisymmetric setting. The statements in that paper
assume $C^2$ boundaries, $a\in W^{3/2,2}(\partial\Omega)$, and
$F\in W^{1,2}(\Omega)$. The planar result imposes no reflection
symmetry. In the spatial result the domain and data are axisymmetric,
and swirl is allowed. These conclusions build on the Euler pressure
and Bernoulli analysis developed in \cite{KPR2013,KPR2015axis}.

A further geometric development is the work of Korobkov, Lyu and Weng
\cite[Theorem~1.1]{KLW2022}. They prove existence for helically
invariant data on the periodic helical domains considered there,
\[
 \Omega\subset\mathbb R^2\times\mathbb T_\sigma,
 \qquad \mathbb T_\sigma=\mathbb R/\sigma\mathbb Z,
 \qquad \sigma>0,
\]
under $C^2$ boundary regularity, $a\in W^{3/2,2}(\partial\Omega)$,
$F\in L^2(\Omega)$, and zero total flux. The axial periodicity is
part of this geometric setting. Their proof combines Leray's
contradiction argument with total-head level sets and the coarea
formula. As explained in \cite[pp.~3976--3977]{KLW2022}, the
stream-function Bernoulli property used in the planar and axisymmetric
analysis can fail in the helical setting, and their argument avoids
that property. The helical theorem concerns a periodic symmetry class
distinct from the bounded planar and axisymmetric classes considered
here.

The monograph of Korobkov, Pileckas and Russo
\cite[Theorems~5.1.1 and~6.1.1]{KPRbook} presents the planar and
axisymmetric general-flux results with $F\in L^2(\Omega)$, while
retaining $C^2$ boundaries and $H^{3/2}$ boundary velocities. These
$L^2$-force formulations provide the benchmark for the force-regularity
comparison in the present work.
The monograph also gives a simplified planar proof in
Sections~5.2--5.3 and retains the Euler-geometric proof for its
axisymmetric extension. Its Remark~6.4.9 places the helical result
\cite{KLW2022} among the subsequent developments of the theory.

Lipschitz boundaries and $H^{1/2}$ boundary velocities already occur
in the existence theory under additional conditions. In
\cite[Theorems~3.3.2, 3.5.1 and~4.1.5]{KPRbook}, these conditions
are, respectively, zero flux through each boundary component, a
small-flux assumption, and an admissible planar reflection symmetry.
For the last result, the domain and data are symmetric and each
boundary component meets the same symmetry line; see
\cite[Definitions~4.1.1--4.1.2]{KPRbook}. These hypotheses are
different from the sole total-flux condition on a general planar
domain. Kozono and Yanagisawa \cite[Theorem~1.3]{KY2009} also use
the dual of the solenoidal energy space, under a variational condition
on the boundary datum.

The present work addresses lower regularity within the bounded planar
and axisymmetric spatial geometries. Theorem~\ref{thm:main} combines
\[
 \partial\Omega\text{ Lipschitz},\qquad
 a\in H^{1/2}(\partial\Omega),\qquad f\in V',
 \qquad \sum_{j=0}^N F_j=0,
\]
where $V$ is defined in \eqref{eq:1.spaces}. The finite
outer-boundary and hole geometry is retained, as is the axisymmetry
of the spatial domain and data. No smallness condition on the
boundary datum or the component fluxes is imposed. The analytic
issue is the simultaneous treatment of the rough physical boundary,
the energy-class trace, and the passage to energy-dual forces in the
general-flux problem. The choice of $V'$ is dictated by the solenoidal
weak formulation. Proposition~\ref{prop:quotient} identifies this
space with the quotient of $H^{-1}$ by pressure gradients. The
force-regularity assertion is therefore the passage from regular
representatives to energy data, with the quotient formulation making
the dependence of the pressure on the representative explicit. The
proof below is developed for the planar and axisymmetric classes
just specified.

\subsection{Solenoidal forcing and the main theorem}

Throughout the paper we distinguish the intrinsic force $f$ from a
representative $F$ acting on all zero-trace test fields. Set
\begin{equation}\label{eq:1.spaces}
\begin{aligned}
 H&=H_0^1(\Omega;\mathbb R^n),
 &\|\varphi\|_H&=\|\nabla\varphi\|_2,\\
 V&=\{\varphi\in H:\operatorname{div}\varphi=0\},
 &H^{-1}&=H',\qquad V'=(V)',\\
 L_0^2(\Omega)&=\left\{q\in L^2(\Omega):\int_\Omega q=0\right\}.
\end{aligned}
\end{equation}
The norm on $V$ is the restriction of the Dirichlet norm on $H$.
Duality brackets are always used with a test field in the domain of
the indicated functional. In particular, $\langle f,\varphi\rangle$
is not used for general $\varphi\in H$ when only $f\in V'$ is given.

For an orthogonal transformation $R$ preserving $\Omega$, write
$\rho_R\varphi(x)=R\varphi(R^{-1}x)$. In the spatial case the symmetry
group $G$ is the rotations about the $x_3$-axis. For the no-swirl
assertion it is enlarged by reflection in a plane containing that axis.
The hypotheses on $a$ and $f$ are
\begin{equation}\label{eq:1.symmetry}
 \rho_R a=a,
 \qquad \langle f,\rho_R\eta\rangle=\langle f,\eta\rangle
 \quad(\eta\in V,\ R\in G).
\end{equation}
For ordinary axisymmetric vector fields, invariance under the additional
reflection is equivalent to vanishing azimuthal component.

\begin{theorem}[Energy trace and solenoidal-dual force]\label{thm:main}
Let $\Omega$ have the geometry above. Assume either $n=2$, or $n=3$
with $\Omega$ invariant under rotations about the $x_3$-axis and the
data satisfying \eqref{eq:1.symmetry}. For
\[
 \nu>0,\qquad a\in H^{1/2}(\partial\Omega;\mathbb R^n),
 \qquad f\in V',\qquad \sum_{j=0}^N F_j=0,
\]
there is a velocity $u\in H^1(\Omega;\mathbb R^n)$ with
$\operatorname{div}u=0$ and $\operatorname{Tr}u=a$ such that
\begin{equation}\label{eq:1.weak}
 \nu\int_\Omega\nabla u:\nabla\eta
 -\int_\Omega(u\otimes u):\nabla\eta
 =\langle f,\eta\rangle
 \qquad(\eta\in V).
\end{equation}
The spatial velocity can be chosen axisymmetric, and can be chosen
without swirl under the additional reflection invariance.

For every selected extension $F\in H^{-1}$ with $F|_V=f$, the same
velocity has a unique associated mean-zero pressure
$P_F\in L_0^2(\Omega)$ satisfying
\begin{equation}\label{eq:1.full}
 \nu\int_\Omega\nabla u:\nabla\varphi
 -\int_\Omega(u\otimes u):\nabla\varphi
 -\int_\Omega P_F\operatorname{div}\varphi
 =\langle F,\varphi\rangle
 \qquad(\varphi\in H).
\end{equation}
Here uniqueness of $P_F$ refers to a fixed velocity $u$ and a fixed
extension $F$; no uniqueness of the nonlinear velocity problem is asserted.
\end{theorem}

\subsection{An isometric extension and the pressure quotient}
\label{subsec:extension}

\begin{lemma}[Canonical energy extension]\label{lem:extension}
There is a linear isometry $\mathcal E:V'\to H^{-1}$ such that
$(\mathcal Ef)|_V=f$. It commutes with every orthogonal symmetry of
$\Omega$ and is given by
\begin{equation}\label{eq:1.extension}
 \mathcal Ef=-\Delta z_f,\qquad
 z_f\in V,\qquad
 (\nabla z_f,\nabla\eta)=\langle f,\eta\rangle\quad(\eta\in V).
\end{equation}
In particular,
\begin{equation}\label{eq:1.extensionnorm}
 \|\mathcal Ef\|_{H^{-1}}=\|\nabla z_f\|_2=\|f\|_{V'}.
\end{equation}
\end{lemma}
\begin{proof}
Since $V$ is a closed subspace of the Hilbert space $H$, Riesz
representation gives a unique $z_f\in V$ and the second equality in
\eqref{eq:1.extensionnorm}. Define
$\langle\mathcal Ef,\varphi\rangle=(\nabla z_f,\nabla\varphi)$ for
all $\varphi\in H$. Cauchy--Schwarz gives
$\|\mathcal Ef\|_{H^{-1}}\le\|\nabla z_f\|_2$; testing with $z_f$
gives equality when $z_f\ne0$, and the zero case is immediate.
Linearity follows from uniqueness. The action $\rho_R$ is unitary
for the Dirichlet inner product and preserves $V$, so the same
uniqueness shows that $\mathcal E$ commutes with that action on
vectors and duals. Thus \eqref{eq:1.symmetry} passes from $f$ to
$\mathcal Ef$ on all of $H$.
\end{proof}

\begin{proposition}[Representatives and pressures]\label{prop:quotient}
On a bounded connected Lipschitz domain, restriction induces an
isometric identification
\begin{equation}\label{eq:1.quotient}
 V'\cong H^{-1}/\nabla L_0^2(\Omega),\qquad
 \|f\|_{V'}=\min_{F|_V=f}\|F\|_{H^{-1}}.
\end{equation}
For any velocity satisfying \eqref{eq:1.weak}, the pressures in
\eqref{eq:1.full} exist and satisfy
\begin{equation}\label{eq:1.pressureestimate}
 \|P_F\|_2\le C_\Omega
 \bigl(\nu\|\nabla u\|_2+\|u\|_4^2+\|F\|_{H^{-1}}\bigr).
\end{equation}
If $F_2-F_1=\nabla q$ with $q\in L_0^2$, then
\begin{equation}\label{eq:1.pressuregauge}
 P_{F_2}=P_{F_1}+q.
\end{equation}
\end{proposition}
\begin{proof}
Let $\mathcal B:L_0^2\to H$ be a bounded right inverse of divergence
\cite{Galdi}. If $\ell\in H'$ vanishes on $V$, then
\[
 \varphi-\mathcal B(\operatorname{div}\varphi)\in V,
 \qquad
 \ell(\varphi)=\ell(\mathcal B(\operatorname{div}\varphi)).
\]
Riesz representation on $L_0^2$ therefore gives a unique $Q\in L_0^2$
such that $\ell(\varphi)=(Q,\operatorname{div}\varphi)$ and
$\|Q\|_2\le\|\mathcal B\|\|\ell\|_{H'}$.
Equivalently, the kernel of restriction is $\nabla L_0^2$.
Lemma~\ref{lem:extension} makes restriction surjective and attains
the minimum extension norm; every extension has norm at least
$\|f\|_{V'}$. This proves \eqref{eq:1.quotient}.

For the pressure, apply this representation to
\[
 \ell_F(\varphi)=\nu(\nabla u,\nabla\varphi)
 -(u\otimes u,\nabla\varphi)-\langle F,\varphi\rangle.
\]
The embedding $H^1\hookrightarrow L^4$ in dimensions two and three
makes $\ell_F\in H'$, while \eqref{eq:1.weak} makes it vanish on $V$.
This gives \eqref{eq:1.full} and \eqref{eq:1.pressureestimate}.
Subtracting the two equations and using
$\langle\nabla q,\varphi\rangle=-(q,\operatorname{div}\varphi)$
gives \eqref{eq:1.pressuregauge}.
\end{proof}

The velocity equation thus depends only on the force class modulo
gradients. The extension in Lemma~\ref{lem:extension} is a convenient
fixed representative, not an additional assumption on the data.
Equivalently it is defined using the orthogonal projection of $H$
onto $V$ for the Dirichlet inner product; this projection is not being
identified with the usual $L^2$ Helmholtz projection.

\begin{remark}[A prescribed distribution is different from a dual class]
\label{rem:rawforce}
Suppose a specified $F_{\rm raw}\in\mathcal D'(\Omega;\mathbb R^n)$
acts continuously on smooth solenoidal test fields for the $V$ norm.
By the density of compactly supported smooth solenoidal fields in
$V$ \cite{Galdi}, it determines $f\in V'$. The de Rham lemma gives
$F_{\rm raw}=\mathcal Ef+\nabla q$ for a distribution $q$.
The pressure for this specified distribution is then
$P_{\mathcal Ef}+q$, and need not lie in $L^2$.
For example $F_{\rm raw}=\nabla\delta_{x_0}$ induces the zero element
of $V'$, whereas $u=0$ in the zero-boundary problem requires
$P=\delta_{x_0}$ modulo constants. Theorem~\ref{thm:main} asserts
$L^2$ regularity relative to $H^{-1}$ representatives, not relative
to arbitrary distributional representatives.
\end{remark}

\subsection{Proof strategy and organization}
\label{subsec:strategy}

The lower regularity of the boundary and trace requires estimates
that do not rely on global $W^{2,2}$ regularity of the viscous
solutions or on smooth distance-parallel boundary surfaces. We use
transverse Lipschitz collars, keep the pressure representative fixed,
and separate the regular-force existence argument from the final
passage to energy-dual data. These choices preserve the geometric
structure of the total-head contradiction while addressing the
boundary and force errors at the energy level.

We use the Euler boundary-pressure and level-set results of
\cite{KPR2013,KPR2015axis,KPR2015}, the Sobolev Morse--Sard theorem
\cite{BKK}, and the Lipschitz Stokes regularity estimate
\cite[Theorem 1.1.3]{MW2012}. Their roles are specified in
Sections~\ref{sec:pressure} and \ref{sec:euler}.
For comparison, the helical proof in \cite{KLW2022} supplies a
different use of total-head level sets; the Euler-geometric inputs
below come from the planar and axisymmetric results just cited.
First we establish the following intermediate result on the fixed domain.

\begin{theorem}[Intermediate existence for regular representatives]
\label{thm:regular}
Under the geometric, trace and symmetry hypotheses of
Theorem~\ref{thm:main}, let $F\in H^1(\Omega;\mathbb R^n)$ be a
specific force, invariant under the required symmetries.
There are $u\in H^1$ and $P_F\in L_0^2$ satisfying
\eqref{eq:1.full}, with the prescribed trace and velocity symmetries.
\end{theorem}

Sections~\ref{sec:continuation}--\ref{sec:regular-existence} prove
Theorem~\ref{thm:regular}. In a normalized blow-up equation write
$\mu_k\to0$ for the viscosity and
$\Phi_k=p_k+|u_k|^2/2$ for the total head. With $c_j$ the limiting
Euler pressure on an admissible boundary component, the collar
estimates give
\[
 \lim_{\delta\downarrow0}\limsup_k
 \operatorname*{ess\,sup}_{\mu_k^2\le h\le\delta}
 \|\Phi_k-c_j\|_{L^2(\Gamma_{j,h})}=0.
\]
Section~\ref{sec:energyforce} fixes $F=\mathcal Ef$ once and for all
and chooses smooth representatives $F_k\to F$ in the full
$H^{-1}$ norm. For a hypothetical unbounded solution sequence, a
small nonlinear correction leaves an effective force satisfying
\[
 \mu_k^{-1}\|e_k\|_{L^{3/2}(\Omega)}\to0,
 \qquad \|e_k\|_{L^2(L_i)}\to0
 \quad\text{for every fixed interior separating layer }L_i.
\]
The first limit controls the volume force error; the second controls
the surface force error after coarea integration in the head level.
The pressure convergence at the end is relative to this fixed
representative. Proposition~\ref{prop:quotient} then supplies the
pressure for any other $H^{-1}$ extension of $f$.

\section{Fixed-domain continuation for regular forces}
\label{sec:continuation}

Let $\Omega\subset\mathbb R^n$, $n=2$ or $3$, be bounded and
connected. We use the geometric description
\[
 \Omega=\Omega_0\setminus\bigcup_{j=1}^N\overline{\Omega_j},
 \qquad \overline{\Omega_j}\subset\Omega_0,
 \qquad \partial\Omega_j=\Gamma_j\quad(0\le j\le N),
\]
where the holes have pairwise disjoint closures and all these
boundaries are Lipschitz and connected. The holes and the exterior
of $\Omega_0$ are connected. In dimension three the domain and
data are axisymmetric about the $x_3$-axis; swirl is allowed.
Each complementary component is invariant under the connected
rotation group, since its orbit in a finite set is a singleton.

For Sections~\ref{sec:continuation}--\ref{sec:regular-existence}
fix a specific representative
\[
 \nu>0,\qquad a\in H^{1/2}(\partial\Omega;\mathbb R^n),
 \qquad F\in H^1(\Omega;\mathbb R^n),\qquad \sum_jF_j=0.
\]
The distinction between the force $F$ and the scalar boundary
fluxes $F_j$ will be maintained. We seek
\begin{equation}\label{eq:2.1}
 -\nu\Delta u+(u\cdot\nabla)u+\nabla P=F,
 \qquad \operatorname{div}u=0,
 \qquad \operatorname{Tr}u=a,
\end{equation}
with $u\in H^1$ and $P\in L_0^2$. Equations are distributional and
boundary values are traces.

\subsection{The solenoidal lift and the weak pressure}

Use $M:N=\sum_{i,j}M_{ij}N_{ij}$ and
$(u\otimes w)_{ij}=u_iw_j$, with $(\nabla\varphi)_{ij}=\partial_j\varphi_i$.
Write $E_{\rm tr}$ for a bounded right inverse of the trace and
$\mathcal B:L_0^2\to H$ for a bounded right inverse of divergence.
Gauss--Green and the total-flux condition give
\[
 \int_\Omega\operatorname{div}(E_{\rm tr}a)
 =\langle\operatorname{Tr}(E_{\rm tr}a),n_\Omega\rangle_{\partial\Omega}
 =\sum_jF_j=0.
\]
Consequently $L=E_{\rm tr}a-\mathcal B(\operatorname{div}E_{\rm tr}a)$
satisfies
\begin{equation}\label{eq:2.2}
 \operatorname{Tr}L=a,\qquad \operatorname{div}L=0,
 \qquad \|L\|_{H^1}\le C_\Omega\|a\|_{H^{1/2}(\partial\Omega)}.
\end{equation}
Let $s\in V$ solve
\[
 \nu(\nabla s,\nabla\eta)
 =\langle F,\eta\rangle-\nu(\nabla L,\nabla\eta)
 \qquad(\eta\in V).
\]
Lax--Milgram and testing with $s$ give, with $U=L+s$,
\begin{equation}\label{eq:2.3}
\begin{gathered}
 \|\nabla s\|_2\le\nu^{-1}\|F|_V\|_{V'}+\|\nabla L\|_2,\\
 \operatorname{Tr}U=a,\qquad\operatorname{div}U=0,
 \qquad \nu(\nabla U,\nabla\eta)=\langle F,\eta\rangle
 \quad(\eta\in V).
\end{gathered}
\end{equation}
In particular $\|F|_V\|_{V'}\le\|F\|_{H^{-1}}$; the Stokes lift
only uses the solenoidal part of the force. We record the pressure
reconstruction from Proposition~\ref{prop:quotient} in the form
\begin{equation}\label{eq:2.4}
 \ell|_V=0\quad\Longrightarrow\quad
 \ell(\varphi)=\int_\Omega P\operatorname{div}\varphi,
 \qquad P\in L_0^2,
 \qquad \|P\|_2\le\|\mathcal B\|\|\ell\|_{H'}.
\end{equation}

\subsection{The compact map and its symmetry}

Set
\[
 b(v,w,z)=\sum_{i,j}\int_\Omega v_j\partial_jw_i\,z_i,
 \qquad |b(v,w,z)|\le\|v\|_4\|\nabla w\|_2\|z\|_4.
\]
If $\operatorname{div}v=0$ and either $w$ or $z$ has zero trace,
integration by parts followed by density gives
\begin{equation}\label{eq:2.5}
 b(v,w,z)=-b(v,z,w).
\end{equation}
For $w\in V$, define $T(w)\in V$ by
\begin{equation}\label{eq:2.6}
 \nu(\nabla T(w),\nabla\eta)
 =-b(U+w,U+w,\eta)
 =\int_\Omega(U+w)\otimes(U+w):\nabla\eta.
\end{equation}
For $y_i=U+w_i$,
\[
 y_1\otimes y_1-y_2\otimes y_2
 =(y_1-y_2)\otimes y_1+y_2\otimes(y_1-y_2),
\]
whence
\[
 \nu\|\nabla(T(w_1)-T(w_2))\|_2
 \le(\|y_1\|_4+\|y_2\|_4)\|w_1-w_2\|_4.
\]
Thus $T$ is continuous. Since $H^1$ embeds compactly into $L^4$
in dimensions two and three, it is also compact. The homotopy is
\begin{equation}\label{eq:2.7}
 w=\lambda T(w)\quad\Longleftrightarrow\quad
 \nu(\nabla w,\nabla\eta)=-\lambda b(U+w,U+w,\eta),
 \qquad 0\le\lambda\le1.
\end{equation}
In space use the group $G$ defined in Section~\ref{sec:main}.
Its action preserves $\Omega$, the $V$ inner product and $b$.
The field $U$ is unique in the affine solenoidal class with trace
$a$: the difference of two such Stokes lifts lies in $V$, and
testing its homogeneous equation with that difference gives zero.
Thus $U$ is independent of the possibly nonsymmetric auxiliary
lift $L$. This uniqueness and uniqueness in \eqref{eq:2.6} imply
\[
 \rho_RU=U,\qquad T(\rho_Rw)=\rho_RT(w),
 \qquad T(V^G)\subset V^G,
 \qquad V^G=\bigcap_{R\in G}\ker(I-\rho_R).
\]
For a $G$-invariant residual $\ell_w$ and $w\in V^G$,
\[
 \ell_w(\eta)=\ell_w\left(\int_G\rho_R\eta\,dR\right),
 \qquad \int_G\rho_R\eta\,dR\in V^G,
\]
where Haar measure is normalized. Testing on $V^G$ therefore
implies the equation on all of $V$. Reflection fixes $e_r,e_z$
and reverses $e_\theta$, giving the no-swirl conclusion when
that reflection is imposed. Below $X$ denotes $V$ in the plane
and the selected fixed subspace in space.

\subsection{Energy identity before the Euler limit}

Suppose $\{w\in X:w=\lambda T(w),\ 0\le\lambda\le1\}$ is
unbounded. Choose
\[
 J_k=\|\nabla w_k\|_2\to\infty,
 \qquad v_k=J_k^{-1}w_k,
 \qquad \|\nabla v_k\|_2=1.
\]
Test \eqref{eq:2.7} with $w_k$. Formula~\eqref{eq:2.5} yields
\[
\begin{gathered}
 b(U,w_k,w_k)=b(w_k,w_k,w_k)=0,
 \\ b(U,U,w_k)=-b(U,w_k,U),
 \qquad b(w_k,U,w_k)=-b(w_k,w_k,U).
\end{gathered}
\]
Division by $J_k^2$ gives
\begin{equation}\label{eq:2.8}
 \nu=\lambda_k b(v_k,v_k,U)
       +(\lambda_k/J_k)b(U,v_k,U).
\end{equation}
Extract $v_k\rightharpoonup v$ in $V$, $v_k\to v$ in $L^4$,
and $\lambda_k\to\lambda_0\in[0,1]$. Since
\[
 |b(v_k,v_k,U)-b(v,v,U)|
 \le\|\nabla U\|_2(\|v_k\|_4+\|v\|_4)\|v_k-v\|_4\to0
\]
and $J_k^{-1}|b(U,v_k,U)|\le J_k^{-1}\|U\|_4^2\|\nabla v_k\|_2\to0$,
we have
\begin{equation}\label{eq:2.9}
 0<\nu=\lambda_0 b(v,v,U),\qquad \lambda_0>0,\qquad v\ne0.
\end{equation}

\subsection{Rescaling the velocity, pressure and force}

Put $\widehat u_k=U+w_k$. The functional
\[
 \ell_k(\varphi)=\nu(\nabla\widehat u_k,\nabla\varphi)
 +\lambda_k b(\widehat u_k,\widehat u_k,\varphi)
 -\langle F,\varphi\rangle
\]
vanishes on $V$, so \eqref{eq:2.4} supplies
$\widehat P_k\in L_0^2$ with
$-\nu\Delta\widehat u_k+\lambda_k(\widehat u_k\cdot\nabla)\widehat u_k
+\nabla\widehat P_k=F$.
Discarding finitely many indices, $\lambda_k\ge\lambda_0/2$.
Define
\[
 u_k=J_k^{-1}\widehat u_k,\qquad
 p_k=(\lambda_kJ_k^2)^{-1}\widehat P_k,\qquad
 \mu_k=\frac\nu{\lambda_kJ_k},\qquad
 g_k=\frac F{\lambda_kJ_k^2}.
\]
Then
\begin{equation}\label{eq:2.10}
 -\mu_k\Delta u_k+(u_k\cdot\nabla)u_k+\nabla p_k=g_k,
 \quad\operatorname{div}u_k=0,
 \quad\operatorname{Tr}u_k=\frac{\lambda_k\mu_k}{\nu}a.
\end{equation}
The trace has one power of $\mu_k$ and the force has two:
\[
 u_k=v_k+J_k^{-1}U,\qquad
 \mu_k^{-1}\int_{\Gamma_j}u_k\cdot n_\Omega\,dS
 =\frac{\lambda_k}{\nu}F_j,\qquad
 g_k=\frac{\lambda_k\mu_k^2}{\nu^2}F.
\]
Consequently
\begin{equation}\label{eq:2.11}
 \mu_k\to0,\qquad \sup_k\|u_k\|_{H^1}<\infty,
 \qquad u_k\to v\text{ in }L^4,
 \qquad \|g_k\|_{H^1}\le\nu^{-2}\mu_k^2\|F\|_{H^1}.
\end{equation}

\subsection{Removal of the nonzero trace from the Stokes equation}

Take viscosity one and force zero in \eqref{eq:2.3}. The resulting
Stokes pair $(A,\pi)$, normalized by $\int_\Omega\pi=0$, satisfies
\[
 -\Delta A+\nabla\pi=0,\qquad\operatorname{div}A=0,
 \qquad\operatorname{Tr}A=a,
\]
\[
 \|A\|_{H^1}+\|\pi\|_2\le C_\Omega\|a\|_{H^{1/2}},
 \qquad\Delta\pi=0\quad\text{in }\mathcal D'(\Omega).
\]
Set
\begin{equation}\label{eq:2.12}
 z_k=u_k-J_k^{-1}A,\qquad
 r_k=(\mu_k/J_k)\pi=(\lambda_k\mu_k^2/\nu)\pi,
 \qquad q_k=p_k-r_k.
\end{equation}
Both $q_k$ and $r_k$ have mean zero. Subtraction gives
\[
 -\Delta(\mu_kz_k)+\nabla q_k
 =-\mu_k\Delta u_k+\nabla p_k
   +(\mu_k/J_k)(\Delta A-\nabla\pi),
\]
and $\|\mu_kz_k\|_{H^1}\le\mu_k(\|u_k\|_{H^1}
+J_k^{-1}\|A\|_{H^1})\to0$. Hence
\begin{equation}\label{eq:2.13}
 -\Delta(\mu_kz_k)+\nabla q_k
 =F_k:=g_k-(u_k\cdot\nabla)u_k,
 \qquad z_k\in V,
 \qquad \|r_k\|_2\le C\mu_k^2.
\end{equation}
Here $F_k$ denotes the normalized Stokes right-hand side in
Sections~\ref{sec:pressure}--\ref{sec:regular-existence}; the notation
will be reset explicitly for the unscaled smooth representatives in
Section~\ref{sec:energyforce}.
Since $1/(3/2)=1/6+1/2$,
\[
 \|F_k\|_{3/2}\le\|g_k\|_{3/2}
                  +\|u_k\|_6\|\nabla u_k\|_2\le C.
\]
In the plane the same calculation gives an $L^s$ bound for every
$1<s<2$, using $1/s=1/t+1/2$ with $t=2s/(2-s)$.
In space, on a rotational collar $r\ge r_*>0$, the cylindrical
components have bounded meridional $H^1(r,z)$ norms; angular
derivatives contribute terms bounded by $r_*^{-1}|u_k|$.
The same improved local $L^s$ bound follows. Constants may depend
on the fixed data and, for these local estimates, on $r_*$.

\section{Pressure compactness at a Lipschitz boundary}\label{sec:pressure}

Write $\mathcal A=\{(0,0,z):z\in\mathbb R\}$ in dimension three. A boundary component $\Gamma$ is called \emph{admissible} if
\[
 n=2,\qquad\text{or}\qquad n=3\ \text{and}\ \rho_\Gamma:=\operatorname{dist}(\Gamma,\mathcal A)>0.
\]
The normalized spatial velocities and pressures are axisymmetric. Constants may depend on $\Omega$ and, for a spatial admissible component, on $\rho_\Gamma$. The collar of Lemma~\ref{lemma:4.1} is denoted by
\[
 X_h(s)=X(s,h),\qquad \Gamma_h=X_h(\Gamma),\qquad
 C_\delta=\{X_h(s):s\in\Gamma,\ 0<h<\delta\}.
\]
Its Jacobians $J_h$ and $J_X$ satisfy
\begin{equation}\label{eq:3.1}
\begin{gathered}
 c\le J_h,J_X\le C,\qquad |\partial_hX|\le C,\\
 |X_h(s)-s|\le Ch,\qquad
 ch\le\operatorname{dist}(X_h(s),\partial\Omega)\le Ch,\\
 c\|Z\circ X_h\|_{L^2(\Gamma)}^2
 \le\|Z\|_{L^2(\Gamma_h)}^2
 \le C\|Z\circ X_h\|_{L^2(\Gamma)}^2.
\end{gathered}
\end{equation}
Suprema over $h$ mean essential suprema unless a continuous interior representative is being used. In dimension three the collar is rotation-equivariant.

\subsection{The homogeneous Stokes regularity estimate}

Fix an aperture $a_*>1$ sufficiently large for the collars and set
\[
 K(s)=\{x\in\Omega:|x-s|<a_*\operatorname{dist}(x,\partial\Omega)\},
 \qquad NZ(s)=\sup_{x\in K(s)}|Z(x)|.
\]
The estimate used from Mitrea--Wright \cite[Theorem~1.1.3]{MW2012} is the following $L^2$ regularity estimate. On a bounded connected Lipschitz domain, if
\[
 b\in W^{1,2}(\partial\Omega;\mathbb R^n),\qquad
 \int_{\Gamma_j}b\cdot n_\Omega\,dS=0\quad(0\le j\le N),
\]
the homogeneous Stokes solution
\[
 -\Delta B+\nabla\Sigma=0,\qquad \operatorname{div}B=0,\qquad
 \operatorname{Tr}B=b,\qquad \int_\Omega\Sigma\,dx=0
\]
has almost everywhere non-tangential pressure limits and satisfies
\begin{equation}\label{eq:3.2}
 \|N(\nabla B)\|_{L^2(\partial\Omega)}
 +\|N\Sigma\|_{L^2(\partial\Omega)}
 \le C\|b\|_{W^{1,2}(\partial\Omega)}.
\end{equation}
Componentwise zero flux is a sufficient condition used for the auxiliary data only. The original datum $a$ still has arbitrary compatible component fluxes. We use the full $W^{1,2}$ norm.

Here is the identification of this solution with the energy solution used below. The layer-potential construction in \cite[Chapter~4 and the proof of the regularity theorem in Section~9.2]{MW2012} gives a pressure representative with almost everywhere non-tangential limits and the maximal estimate above. For either $Z=\nabla B$ or $Z=\Sigma$ in this gauge, integration along the collars of Lemma~\ref{lemma:4.1} gives
\[
 \int_{C_{h_0}}|Z|^2\,dx
 \le Ch_0\int_{\partial\Omega}|NZ|^2\,dS.
\]
On the compact interior remainder, each point belongs to the cones over a boundary patch of uniformly positive surface measure: choose a closest boundary point and use the lower area bound for Lipschitz graphs. Thus
\[
 \|Z\|_{L^2(\Omega)}\le C_\Omega\|NZ\|_{L^2(\partial\Omega)}.
\]
Subtracting the volume mean from $\Sigma$ preserves the estimate up to a domain constant. Integrating $\nabla B$ along the same paths from its boundary value $b$ gives $B\in L^2(\Omega)$, and the Sobolev trace is $b$. If two finite-energy pairs have this trace, their velocity difference $d$ belongs to $V$ and solves the homogeneous weak equation. Testing by $d$ gives $\|\nabla d\|_2^2=0$; their pressure gradients then agree. Connectedness and the mean-zero convention identify their pressures.

For pressures with arbitrary gauges put $\langle\Sigma\rangle=|\Omega|^{-1}\int_\Omega\Sigma$. Applying \eqref{eq:3.2} to a difference gives
\begin{equation}\label{eq:3.3}
\begin{split}
 \|N(\Sigma_1-\Sigma_2)\|_{L^2(\partial\Omega)}
 &\le C\|b_1-b_2\|_{W^{1,2}(\partial\Omega)}\\
 &\quad+\mathcal H^{n-1}(\partial\Omega)^{1/2}
       |\langle\Sigma_1-\Sigma_2\rangle|.
\end{split}
\end{equation}
Thus both the boundary data and the pressure constants must converge. For all $0<h<h_0$, \eqref{eq:3.1} also gives
\begin{equation}\label{eq:3.4}
 X_h(s)\in K(s),\qquad
 \|Z\circ X_h\|_{L^2(\Gamma)}\le\|NZ\|_{L^2(\Gamma)}.
\end{equation}

\subsection{Strong bulk convergence without boundary derivative estimates}

The normalized sequence has the properties
\begin{equation}\label{eq:3.5}
\begin{gathered}
 \sup_k\|u_k\|_{H^1(\Omega)}<\infty,\qquad
 \operatorname{div}u_k=0,\qquad
 \operatorname{Tr}u_k\longrightarrow0\quad\text{in }H^{1/2},\\
 g_k\longrightarrow0\quad\text{in }H^{-1}(\Omega),\qquad
 \sup_k\|g_k\|_{L^{3/2}(\Omega)}<\infty.
\end{gathered}
\end{equation}
Define $F_k=g_k-(u_k\cdot\nabla)u_k$. H\"older's inequality gives
\begin{equation}\label{eq:3.6}
 \|F_k\|_{3/2}
 \le\|g_k\|_{3/2}+\|u_k\|_6\|\nabla u_k\|_2
 \le\|g_k\|_{3/2}+C\|u_k\|_{H^1}^2\le C.
\end{equation}
The pressure separation in \eqref{eq:2.12} yields
\begin{equation}\label{eq:3.7}
\begin{gathered}
 Y_k=\mu_kz_k\in V,\qquad -\Delta Y_k+\nabla q_k=F_k,\\
 \int_\Omega q_k\,dx=0,\qquad Y_k\longrightarrow0\quad\text{in }H^1(\Omega).
\end{gathered}
\end{equation}

\begin{lemma}[Bulk pressure convergence]\label{lemma:3.1}
After extraction,
\[
 u_k\rightharpoonup v\ \text{in }H^1,\qquad
 u_k\longrightarrow v\ \text{in }L^4,\qquad
 q_k\longrightarrow p\ \text{in }L^2,
\]
where
\[
 v\in V,\qquad \int_\Omega p=0,\qquad
 \nabla p=-(v\cdot\nabla)v,\qquad p\in W^{1,3/2}(\Omega).
\]
If $n=2$, then $p\in W^{1,t}(\Omega)$ for every $1<t<2$.
\end{lemma}

\begin{proof}
Rellich compactness and the trace theorem imply
\[
 u_k\rightharpoonup v\ \text{in }H^1,\qquad
 u_k\longrightarrow v\ \text{in }L^4,\qquad
 \operatorname{div}v=0,\qquad\operatorname{Tr}v=0.
\]
For $\varphi\in H_0^1(\Omega;\mathbb R^n)$,
\[
 \langle(u_k\cdot\nabla)u_k,\varphi\rangle
 =-\int_\Omega(u_k\otimes u_k):\nabla\varphi\,dx.
\]
Consequently, with $F=-\operatorname{div}(v\otimes v)$,
\begin{equation}\label{eq:3.8}
\begin{split}
 \|F_k-F\|_{H^{-1}}
 &\le\|g_k\|_{H^{-1}}+\|u_k\otimes u_k-v\otimes v\|_2\\
 &\le\|g_k\|_{H^{-1}}+(\|u_k\|_4+\|v\|_4)\|u_k-v\|_4
 \longrightarrow0.
\end{split}
\end{equation}
Here is the continuity estimate for the pressure. If
\[
 -\Delta Z+\nabla\rho=H,\qquad Z\in V,\qquad\rho\in L_0^2(\Omega),
\]
testing with $Z$ gives
\[
 \|\nabla Z\|_2^2=\langle H,Z\rangle
 \le\|H\|_{H^{-1}}\|\nabla Z\|_2.
\]
The divergence inf-sup inequality yields, for a fixed $\beta_\Omega>0$,
\[
\begin{split}
 \beta_\Omega\|\rho\|_2
 &\le\sup_{\varphi\in H_0^1\setminus\{0\}}
 \frac{|\int_\Omega\rho\operatorname{div}\varphi|}{\|\nabla\varphi\|_2}\\
 &\le\|H\|_{H^{-1}}+\|\nabla Z\|_2
 \le C_\Omega\|H\|_{H^{-1}}.
\end{split}
\]
Applied to $(Y_k-Y_l,q_k-q_l,F_k-F_l)$, these inequalities show
\begin{equation}\label{eq:3.9}
 \|Y_k-Y_l\|_{H^1}+\|q_k-q_l\|_2
 \le C_\Omega\|F_k-F_l\|_{H^{-1}}\longrightarrow0.
\end{equation}
Since $Y_k\to0$, passing to the distributional limit gives
\[
 \nabla p=F=-(v\cdot\nabla)v,\qquad
 \|\nabla p\|_{3/2}\le\|v\|_6\|\nabla v\|_2.
\]
Together with $p\in L^2(\Omega)\subset L^{3/2}(\Omega)$, this proves $p\in W^{1,3/2}$. In two dimensions, for $1<t<2$,
\[
 \frac1t=\frac{2-t}{2t}+\frac12,\qquad
 \|\nabla p\|_t\le\|v\|_{2t/(2-t)}\|\nabla v\|_2<\infty.
\]
\end{proof}

\subsection{Whole-space potentials and the homogeneous remainder}

Let $\widetilde F_k$ be the extension of $F_k$ by zero. The Stokes fundamental velocity tensor and pressure vector can be chosen as
\[
\begin{array}{lll}
 n=3:&\displaystyle E_{ij}(x)=\frac1{8\pi}
       \left(\frac{\delta_{ij}}{|x|}+\frac{x_ix_j}{|x|^3}\right),
      &\displaystyle P_j(x)=\frac{x_j}{4\pi|x|^3},\\[6pt]
 n=2:&\displaystyle E_{ij}(x)=\frac1{4\pi}
       \left(-\delta_{ij}\log|x|+\frac{x_ix_j}{|x|^2}\right),
      &\displaystyle P_j(x)=\frac{x_j}{2\pi|x|^2}.
\end{array}
\]
They satisfy
\[
 -\Delta E_{ij}+\partial_iP_j=\delta_{ij}\delta_0,
 \qquad \sum_i\partial_iE_{ij}=0.
\]
Put
\begin{equation}\label{eq:3.10}
 W_{k,i}=\sum_jE_{ij}*\widetilde F_{k,j},\qquad
 Q_k=\sum_jP_j*\widetilde F_{k,j}.
\end{equation}
Then
\[
 -\Delta W_k+\nabla Q_k=\widetilde F_k,\qquad
 \operatorname{div}W_k=0\quad\text{in }\mathbb R^n.
\]
For $1<s<\infty$, Calder\'on--Zygmund estimates give
\[
 \|D^2W_k\|_{L^s(\mathbb R^n)}+\|\nabla Q_k\|_{L^s(\mathbb R^n)}
 \le C_s\|\widetilde F_k\|_{L^s(\mathbb R^n)}.
\]
The distributional kernel derivatives include the usual bounded identity terms. For a fixed ball $B$ with $\Omega\Subset B$, choose $R$ so that $x-y\in B_R$ for $x\in B$, $y\in\Omega$. Then
\[
 E,\nabla E,P\in L^1(B_R),
\]
and Young's inequality controls all lower derivatives on $B$. In particular,
\begin{equation}\label{eq:3.11}
 \|W_k\|_{W^{2,3/2}(B)}+\|Q_k\|_{W^{1,3/2}(B)}
 \le C_{B,\Omega}\|F_k\|_{L^{3/2}(\Omega)}\le C_{B,\Omega}.
\end{equation}
Only bounded-set norms occur; the logarithmic two-dimensional kernel requires no cancellation of the total force here.

Since
\[
 \frac1{3/2}-\frac1n<\frac12\qquad(n=2,3),
\]
the compact embeddings $W^{2,3/2}(B)\Subset W^{1,2}(B)$ and $W^{1,3/2}(B)\Subset L^2(B)$ imply, on one subsequence,
\begin{equation}\label{eq:3.12}
\begin{gathered}
 W_k\rightharpoonup W\quad\text{in }W^{2,3/2}(B),\qquad
 W_k\longrightarrow W\quad\text{in }H^1(B),\\
 Q_k\rightharpoonup Q\quad\text{in }W^{1,3/2}(B),\qquad
 Q_k\longrightarrow Q\quad\text{in }L^2(B).
\end{gathered}
\end{equation}
For a rotation $R$, the identities
\[
 E(Rx)=RE(x)R^T,\qquad P(Rx)=RP(x)
\]
show that axisymmetric $F_k$ produce $W_k(Rx)=RW_k(x)$ and $Q_k(Rx)=Q_k(x)$.

Define
\begin{equation}\label{eq:3.13}
 b_k=Y_k-W_k|_\Omega,\qquad\sigma_k=q_k-Q_k|_\Omega.
\end{equation}
Subtracting the equations gives
\begin{equation}\label{eq:3.14}
 -\Delta b_k+\nabla\sigma_k=0,\qquad
 \operatorname{div}b_k=0,\qquad\operatorname{Tr}b_k=-\operatorname{Tr}W_k.
\end{equation}
Moreover,
\begin{equation}\label{eq:3.15}
 b_k\longrightarrow b=-W|_\Omega\quad\text{in }H^1,\qquad
 \sigma_k\longrightarrow\sigma=p-Q|_\Omega\quad\text{in }L^2.
\end{equation}
For $K\Subset K_1\Subset\Omega$ and each integer $m\ge0$, homogeneous interior Stokes estimates and Sobolev embedding give
\begin{equation}\label{eq:3.16}
\begin{split}
 &\|b_k-b\|_{C^m(K)}+\|\sigma_k-\sigma\|_{C^m(K)}\\
 &\qquad\le C_{m,K,K_1}
   \bigl(\|b_k-b\|_{H^1(K_1)}+\|\sigma_k-\sigma\|_{L^2(K_1)}\bigr)
 \longrightarrow0.
\end{split}
\end{equation}

\begin{remark}[Interior derivative bounds]\label{remark:3.2}\label{lemma:3.2}
For $K\Subset K_1\Subset\Omega$, the inhomogeneous interior estimate is
\[
\begin{split}
 \|Y_k\|_{W^{2,3/2}(K)}+\|q_k\|_{W^{1,3/2}(K)}
 &\le C_{K,K_1}
 \bigl(\|F_k\|_{3/2}+\|Y_k\|_{H^1}+\|q_k\|_2\bigr)\\
 &\le C_{K,K_1}.
\end{split}
\]
The equation has unit viscosity, so the constant contains no negative power of $\mu_k$. Harmonic interior estimates also give
\[
 \|r_k\|_{C^m(K)}\le C_{m,K}\|r_k\|_{L^2(\Omega)}=O(\mu_k^2).
\]
For $\Phi_k=p_k+\tfrac12|u_k|^2$ and $\Phi=p+\tfrac12|v|^2$,
\[
\begin{gathered}
 \|\Phi_k-\Phi\|_2
 \le\|q_k-p\|_2+\|r_k\|_2
 +\tfrac12(\|u_k\|_4+\|v\|_4)\|u_k-v\|_4\longrightarrow0,\\
 \|\nabla|u_k|^2\|_{3/2}\le2\|u_k\|_6\|\nabla u_k\|_2\le C.
\end{gathered}
\]
Thus $\sup_k\|\Phi_k\|_{W^{1,3/2}(K)}<\infty$. On spatial meridional compact sets separated from the axis the same estimate holds in ordinary two-dimensional Sobolev norms.
\end{remark}

\subsection{Compact traces of the whole-space potentials}

\begin{lemma}\label{lemma:3.3}
In the plane, after extraction,
\[
 \operatorname{Tr}W_k\longrightarrow\operatorname{Tr}W
 \quad\text{in }W^{1,2}(\partial\Omega),\qquad
 \operatorname{Tr}Q_k\longrightarrow\operatorname{Tr}Q
 \quad\text{in }L^2(\partial\Omega).
\]
In the spatial axisymmetric case these convergences hold on every admissible component $\Gamma$.
\end{lemma}

\begin{proof}
In dimension two, for any $Z_k$ equal to a component of $\nabla W_k$ or to $Q_k$, the trace estimate on a Lipschitz curve is
\[
 \|\operatorname{Tr}Z_k\|_{W^{1/3,3/2}(\partial\Omega)}
 \le C\|Z_k\|_{W^{1,3/2}(B)}\le C.
\]
The trace has dimension one, and its critical exponent is
\[
 \frac1{q_*}=\frac23-\frac13=\frac13,\qquad q_*=3.
\]
Consequently, after extraction,
\begin{equation}\label{eq:3.17}
 W^{1/3,3/2}(\partial\Omega)\Subset L^2(\partial\Omega),\qquad
 \begin{cases}
  \operatorname{Tr}\nabla W_k\longrightarrow\operatorname{Tr}\nabla W&\text{in }L^2,\\
  \operatorname{Tr}Q_k\longrightarrow\operatorname{Tr}Q&\text{in }L^2.
 \end{cases}
\end{equation}
The limits are identified using weak continuity of the trace map and \eqref{eq:3.12}. In addition,
\[
 \|\operatorname{Tr}(W_k-W)\|_{L^2(\partial\Omega)}
 \le C\|W_k-W\|_{H^1(B)}\longrightarrow0.
\]
Let $\gamma:I\to\partial\Omega$ be a Lipschitz chart. Smooth ambient approximation gives the weak tangential chain rule
\[
 \partial_s(\operatorname{Tr}W_k\circ\gamma)
   =(\operatorname{Tr}\nabla W_k\circ\gamma)\gamma'
 \quad\text{a.e. on }I.
\]
Indeed, for ambient smooth approximants $W_k^{(m)}$,
\[
 W_k^{(m)}\longrightarrow W_k\ \text{in }W^{2,3/2}(B)
 \quad\Longrightarrow\quad
 \begin{cases}
 \operatorname{Tr}W_k^{(m)}\longrightarrow\operatorname{Tr}W_k
     &\text{in }L^2(\gamma(I)),\\
 \operatorname{Tr}\nabla W_k^{(m)}\longrightarrow\operatorname{Tr}\nabla W_k
     &\text{in }L^2(\gamma(I)).
 \end{cases}
\]
The chain rule for $W_k^{(m)}$ therefore passes to distributions and to $L^2(I)$. In particular,
\[
\begin{split}
 \|\partial_s[\operatorname{Tr}(W_k-W)\circ\gamma]\|_{L^2(I)}
 &\le\|\gamma'\|_\infty
 \|(\operatorname{Tr}\nabla(W_k-W))\circ\gamma\|_{L^2(I)}\\
 &\longrightarrow0.
\end{split}
\]
Summing over finitely many charts proves the planar assertion.

For $n=3$, fix a rotational neighborhood $U$ of $\Gamma$ with $0<\rho\le r\le R$. Write
\[
 W_k=w_{k,r}(r,z)e_r(\theta)+w_{k,\theta}(r,z)e_\theta(\theta)
       +w_{k,z}(r,z)e_z,\qquad Q_k=Q_k(r,z).
\]
For the generating neighborhood $U'$,
\[
\begin{gathered}
 \sum_{\alpha=r,\theta,z}\|w_{k,\alpha}\|_{W^{2,3/2}(U')}
       +\|Q_k\|_{W^{1,3/2}(U')}\le C_{\rho,R},\\
 \|Z\|_{L^s(U)}^s
 =2\pi\int_{U'}|Z(r,z)|^s r\,dr\,dz
 \asymp_{\rho,R}\|Z\|_{L^s(U')}^s.
\end{gathered}
\]
Cartesian and cylindrical derivatives up to order two differ by bounded coefficient combinations of derivatives of lower or equal order; the coefficients contain only bounded powers of $r^{-1}$. Apply the one-dimensional trace compactness in \eqref{eq:3.17} to these scalar components on the generating curve $\gamma(s)=(r(s),z(s))$. With
\[
 \Xi(s,\theta)=(r(s)\cos\theta,r(s)\sin\theta,z(s)),\qquad
 dS=r(s)|\gamma'(s)|\,ds\,d\theta,
\]
the angular derivatives are
\[
 \partial_\theta(W_k\circ\Xi)=w_{k,r}e_\theta-w_{k,\theta}e_r,
 \qquad\partial_\theta(Q_k\circ\Xi)=0.
\]
The meridional derivatives converge by the preceding chain rule; the angular derivatives converge by the $L^2$ convergence of the component values. Since $\rho\le r(s)\le R$, these statements give exactly the asserted physical $W^{1,2}(\Gamma)$ and $L^2(\Gamma)$ convergences. Only traces of first-order Sobolev functions have been used.
\end{proof}

\subsection{Localization at a whole off-axis component}

\begin{lemma}[Compact homogeneous pressure in a collar]\label{lemma:3.4}
For each admissible $\Gamma$, one subsequence and a fixed $h_0>0$ satisfy
\begin{equation}\label{eq:3.18}
 \sup_{0<h<h_0}\|(\sigma_k-\sigma)\circ X_h\|_{L^2(\Gamma)}
 \longrightarrow0.
\end{equation}
The pressure $\sigma$ has an $L^2$ non-tangential trace near $\Gamma$.
\end{lemma}

\begin{proof}
\emph{Planar case.} Each $\Gamma_j$ bounds a bounded Lipschitz region $O_j$. Since $W_k\in H^1(O_j)$ and $\operatorname{div}W_k=0$ there,
\[
 \int_{\Gamma_j}W_k\cdot n_\Omega\,dS
 =\pm\int_{O_j}\operatorname{div}W_k\,dx=0.
\]
Lemma~\ref{lemma:3.3} and \eqref{eq:3.15} give
\[
 \|\operatorname{Tr}(b_k-b)\|_{W^{1,2}(\partial\Omega)}\longrightarrow0,
 \qquad
 |\langle\sigma_k-\sigma\rangle|
 \le|\Omega|^{-1/2}\|\sigma_k-\sigma\|_2\longrightarrow0.
\]
Now \eqref{eq:3.3} implies $\|N(\sigma_k-\sigma)\|_2\to0$, and \eqref{eq:3.4} proves \eqref{eq:3.18}.

\emph{Spatial localization.} Fix an admissible $\Gamma$ and choose $\chi\in C^\infty(\mathbb R^3)$ satisfying
\[
 \chi=1\ \text{near }\Gamma,\qquad
 \chi=0\ \text{near }(\partial\Omega\setminus\Gamma)\cup\mathcal A,
 \qquad0\le\chi\le1.
\]
After restricting to a bounded ambient neighborhood,
\[
 K:=\operatorname{supp}(\nabla\chi)\cap\Omega\Subset\Omega\setminus\mathcal A.
\]
Let $G_k^0=\nabla\chi\cdot b_k$. The product rule and Gauss--Green give
\begin{equation}\label{eq:3.19}
\begin{split}
 \int_\Omega G_k^0\,dx
 &=\int_\Omega\operatorname{div}(\chi b_k)\,dx
   =\int_\Gamma b_k\cdot n_\Omega\,dS\\
 &=-\int_\Gamma W_k\cdot n_\Omega\,dS=0.
\end{split}
\end{equation}
Moreover, \eqref{eq:3.16} gives
\[
 G_k^0\longrightarrow G^0:=\nabla\chi\cdot b\quad\text{in }C^m(\overline\Omega),
 \qquad\operatorname{supp}G_k^0,\operatorname{supp}G^0\subset K
 \quad(m\ge0).
\]

For clarity, a fixed compact-support divergence correction can be constructed as follows. Cover $K$ by finitely many interior balls and add finite chains of interior balls so that their intersection graph is connected. Write the resulting balls as $D_1,\ldots,D_M\Subset\Omega$ and choose a rooted spanning tree, with root $1$ and parent index $\pi(j)$ for $j\ne1$. Choose
\[
 \sum_{j=1}^M\eta_j=1\ \text{near }K,\qquad
 \eta_j\in C_c^\infty(D_j),\qquad
 \vartheta_j\in C_c^\infty(D_j\cap D_{\pi(j)}),\qquad
 \int\vartheta_j=1\quad(j\ne1).
\]
For $G\in C_c^\infty(\Omega)$ supported in $K$ with $\int G=0$, let
\[
 G_j=\eta_jG,\qquad m_j=\int G_j,\qquad
 M_j=\sum_{l\in\operatorname{desc}(j)}m_l,
\]
where the descendants include $j$ itself. Define
\[
\begin{split}
 H_j&=G_j+\sum_{l:\pi(l)=j}M_l\vartheta_l-M_j\vartheta_j
       \quad(j\ne1),\\
 H_1&=G_1+\sum_{l:\pi(l)=1}M_l\vartheta_l.
\end{split}
\]
Then
\[
 \operatorname{supp}H_j\Subset D_j,\qquad
 \int H_j=0,\qquad \sum_{j=1}^M H_j=G.
\]
On $D_j$ choose $\zeta_j\in C_c^\infty(D_j)$, $\int\zeta_j=1$, supported in a smaller concentric ball. One may take
\[
 (\mathcal B_jH_j)(x)
 =\int_{D_j}H_j(y)(x-y)
   \left(\int_1^\infty\zeta_j(y+t(x-y))t^{n-1}\,dt\right)dy.
\]
The mean-zero identity gives $\operatorname{div}\mathcal B_jH_j=H_j$. Its support is contained in the convex hull of $\operatorname{supp}H_j\cup\operatorname{supp}\zeta_j$, a fixed compact subset of $D_j$. Put $\mathcal BG=\sum_j\mathcal B_jH_j$. The standard Bogovski\u\i\ estimates \cite{Galdi} imply, for $m\ge0$ and $1<s<\infty$,
\begin{equation}\label{eq:3.20}
\begin{gathered}
 \operatorname{div}\mathcal BG=G,\qquad
 \operatorname{supp}\mathcal BG\subset K'\Subset\Omega,\\
 \|\mathcal BG\|_{W^{m+1,s}(\Omega)}
 \le C_{m,s,K}\|G\|_{W^{m,s}(\Omega)}.
\end{gathered}
\end{equation}
The compact set $K'$ is fixed because all supports, balls and auxiliary bump functions are fixed. Taking
\[
 d_k=\mathcal BG_k^0,\qquad d=\mathcal BG^0,
\]
and using Sobolev embedding at arbitrarily high finite orders yields
\begin{equation}\label{eq:3.21}
 \operatorname{div}d_k=G_k^0,\qquad
 d_k\longrightarrow d\quad\text{in }C^m(\overline\Omega),\qquad
 \operatorname{supp}d_k\subset K'\Subset\Omega\quad(m\ge0).
\end{equation}
No symmetry of $d_k$ is needed.

Set $B_k=\chi b_k-d_k$ and $\Sigma_k=\chi\sigma_k$. Then
\begin{equation}\label{eq:3.22}
\begin{aligned}
 -\Delta B_k+\nabla\Sigma_k
 &=\chi(-\Delta b_k+\nabla\sigma_k)
   -2\sum_l(\partial_l\chi)\partial_l b_k\\
 &\quad-(\Delta\chi)b_k+\sigma_k\nabla\chi+\Delta d_k\\
 &=-2\nabla\chi\cdot\nabla b_k-(\Delta\chi)b_k
   +\sigma_k\nabla\chi+\Delta d_k=:G_k,\\
 \operatorname{div}B_k&=\nabla\chi\cdot b_k-\operatorname{div}d_k=0,\\
 \operatorname{Tr}B_k&=-\chi\operatorname{Tr}W_k.
\end{aligned}
\end{equation}
With $K''=K\cup K'\Subset\Omega$,
\[
 \operatorname{supp}G_k\subset K'',\qquad
 G_k\longrightarrow G\quad\text{in }C^m(\overline\Omega)\quad(m\ge0),
\]
by \eqref{eq:3.16} and \eqref{eq:3.21}. The boundary data satisfy
\[
 -\chi\operatorname{Tr}W_k
 =\begin{cases}
   -\operatorname{Tr}W_k&\text{on }\Gamma,\\
   0&\text{on }\partial\Omega\setminus\Gamma,
  \end{cases}
 \qquad
 \|\chi\operatorname{Tr}(W_k-W)\|_{W^{1,2}(\partial\Omega)}\longrightarrow0.
\]
Extend $G_k$ by zero and define
\[
 V_k=E*G_k,\qquad T_k=P*G_k.
\]
For a fixed ambient ball $B$, the potential estimates at arbitrary finite Sobolev orders, followed by embedding, give
\[
 V_k\longrightarrow V,\quad T_k\longrightarrow T
 \quad\text{in }C^m(\overline B)\quad(m\ge0).
\]
For example, differentiating the compactly supported sources and using local integrability of the kernels gives
\[
 \|V_k-V\|_{C^m(\overline B)}+\|T_k-T\|_{C^m(\overline B)}
 \le C_{m,B,K''}\|G_k-G\|_{C^m(\mathbb R^3)}\longrightarrow0.
\]
Now
\[
 \widehat B_k=B_k-V_k,\qquad\widehat\Sigma_k=\Sigma_k-T_k
\]
solve homogeneous Stokes with trace
\[
 \widehat b_k=-\chi\operatorname{Tr}W_k-\operatorname{Tr}V_k
 \longrightarrow-\chi\operatorname{Tr}W-\operatorname{Tr}V=:\widehat b
 \quad\text{in }W^{1,2}(\partial\Omega).
\]
Every component flux vanishes:
\[
 \int_{\Gamma_j}\widehat b_k\cdot n_\Omega\,dS
 =-\chi|_{\Gamma_j}\int_{\Gamma_j}W_k\cdot n_\Omega\,dS
   -\int_{\Gamma_j}V_k\cdot n_\Omega\,dS=0.
\]
Here both $W_k$ and $V_k$ are globally solenoidal, so Gauss--Green applies in the bounded region enclosed by each $\Gamma_j$. Furthermore,
\[
 \widehat\Sigma_k\longrightarrow\chi\sigma-T\quad\text{in }L^2(\Omega),
 \qquad\langle\widehat\Sigma_k\rangle\longrightarrow\langle\chi\sigma-T\rangle.
\]
Apply \eqref{eq:3.3} to conclude
\[
 \|N[\widehat\Sigma_k-(\chi\sigma-T)]\|_2\longrightarrow0.
\]
Since $T_k\to T$ uniformly on $\Omega$,
\[
\begin{split}
 \|N(\Sigma_k-\chi\sigma)\|_2
 &\le\|N[\widehat\Sigma_k-(\chi\sigma-T)]\|_2\\
 &\quad+|\partial\Omega|^{1/2}\|T_k-T\|_\infty\longrightarrow0.
\end{split}
\]
Choose $h_0$ so that $\chi=1$ on $C_{h_0}$. Then
\[
 \sup_{0<h<h_0}\|(\sigma_k-\sigma)\circ X_h\|_{L^2(\Gamma)}
 \le\|N(\Sigma_k-\chi\sigma)\|_{L^2(\Gamma)}\longrightarrow0.
\]
The same homogeneous estimate and the smooth potential $T$ give the asserted non-tangential trace of $\sigma$. There are finitely many admissible components, so finitely many extractions suffice.
\end{proof}

\subsection{A quantitative collar inequality and trace identification}

All Sobolev norms in this subsection are on a fixed two-dimensional collar: the physical collar if $n=2$, its meridional representative if $n=3$. Denote it by $C'_{h_0}$. In the spatial case the physical and meridional integrals are equivalent because $0<\rho\le r\le R$.

For $Z\in W^{1,3/2}(C'_{h_0})$,
\begin{equation}\label{eq:3.23}
 \|Z\|_6\le C\|Z\|_{W^{1,3/2}},\qquad
 \|Z\nabla Z\|_{6/5}
 \le\|Z\|_6\|\nabla Z\|_{3/2}
 \le C\|Z\|_{W^{1,3/2}}^2.
\end{equation}
For a smooth $Z$, the fundamental theorem of calculus in the fixed measure on $\Gamma$ gives, for $0<h\le\delta$,
\[
\begin{split}
 &\int_\Gamma|Z(X_h(s))|^2\,dS(s)
       -\int_\Gamma|\operatorname{Tr}Z(s)|^2\,dS(s)\\
 &=2\int_0^h\int_\Gamma
      Z(X_t(s))\nabla Z(X_t(s))\cdot\partial_tX_t(s)\,dS(s)\,dt\\
 &\le C\int_{C'_\delta}|Z||\nabla Z|\,dr\,dz
 \le C|C'_\delta|^{1/6}
             \|Z\nabla Z\|_{L^{6/5}(C'_{h_0})}.
\end{split}
\]
In the planar case the variables in the last integral are the usual Cartesian variables. Since $|C'_\delta|\le C\delta$,
\begin{equation}\label{eq:3.24}
 \sup_{0<h\le\delta}\|Z\circ X_h\|_{L^2(\Gamma)}^2
 \le C\|\operatorname{Tr}Z\|_{L^2(\Gamma)}^2
       +C\delta^{1/6}\|Z\|_{W^{1,3/2}(C'_{h_0})}^2.
\end{equation}
For general $Z$, let $Z^{(m)}\to Z$ in $W^{1,3/2}$ by smooth extension and approximation. The trace map into $L^2(\Gamma)$ is continuous, and \eqref{eq:3.24} applied to $Z^{(m)}-Z^{(l)}$ proves convergence of the flow restrictions in $L^\infty((0,h_0);L^2(\Gamma))$. This defines their representatives and passes the inequality to $Z$.

The same argument identifies the boundary value. For fixed $m$,
\[
\begin{split}
 \|Z\circ X_h-\operatorname{Tr}Z\|_2
 &\le\|(Z-Z^{(m)})\circ X_h\|_2
      +\|Z^{(m)}\circ X_h-\operatorname{Tr}Z^{(m)}\|_2\\
 &\quad+\|\operatorname{Tr}(Z^{(m)}-Z)\|_2.
\end{split}
\]
First let $h\downarrow0$; the middle term vanishes. The other terms are bounded by $C\|Z-Z^{(m)}\|_{W^{1,3/2}}$, uniformly for $0<h<h_0$. Letting $m\to\infty$ proves
\begin{equation}\label{eq:3.25}
 Z\circ X_h\longrightarrow\operatorname{Tr}Z
 \quad\text{in }L^2(\Gamma)\quad(h\downarrow0).
\end{equation}

\begin{lemma}[Pressure compactness at admissible components]\label{lemma:3.5}
Under \eqref{eq:3.5}--\eqref{eq:3.7}, suppose the Euler pressure has constant trace $c_\Gamma$ on an admissible component $\Gamma$. Then
\begin{equation}\label{eq:3.26}
 \lim_{\delta\downarrow0}\limsup_{k\to\infty}
 \sup_{0<h\le\delta}\|q_k-c_\Gamma\|_{L^2(\Gamma_h)}=0.
\end{equation}
One extracted subsequence works for all admissible components.
\end{lemma}

\begin{proof}
From Lemma~\ref{lemma:3.1} and the potential estimates,
\[
 p,Q,\sigma=p-Q\in W^{1,3/2}(C'_{h_0}),\qquad
 \sup_k\|Q_k-Q\|_{W^{1,3/2}(C'_{h_0})}<\infty.
\]
The trace of $\sigma$ has two equivalent descriptions. Formula~\eqref{eq:3.25} gives its Sobolev trace, while Lemma~\ref{lemma:3.4} gives a non-tangential trace $\sigma^{\mathrm{nt}}$. Along the transverse paths,
\[
 \sigma(X_h(s))\longrightarrow\sigma^{\mathrm{nt}}(s)\quad\text{a.e.},
 \qquad |\sigma(X_h(s))|\le N_\Gamma\sigma(s)\in L^2(\Gamma).
\]
Here $N_\Gamma$ may be restricted to the fixed neighborhood on which the localized regularity estimate was obtained. Dominated convergence and uniqueness of an $L^2$ limit imply
\[
 \operatorname{Tr}\sigma=\sigma^{\mathrm{nt}},\qquad
 \operatorname{Tr}Q+\sigma^{\mathrm{nt}}
       =\operatorname{Tr}p=c_\Gamma.
\]
For an explicit estimate set
\[
 e_k=\|\operatorname{Tr}(Q_k-Q)\|_{L^2(\Gamma)},\qquad
 n_k=\sup_{0<h<h_0}\|(\sigma_k-\sigma)\circ X_h\|_{L^2(\Gamma)}.
\]
The preceding lemmas give $e_k\to0$ and $n_k\to0$. Decompose
\[
 q_k-c_\Gamma=(Q_k-Q)+(\sigma_k-\sigma)+(p-c_\Gamma).
\]
Apply \eqref{eq:3.24} to the first and third terms; the trace of the third is zero. With a constant $M_\Gamma<\infty$ independent of $k,h,\delta$,
\begin{equation}\label{eq:3.27}
\begin{gathered}
 \sup_{0<h\le\delta}\|q_k-c_\Gamma\|_{L^2(\Gamma_h)}^2
 \le C(e_k^2+n_k^2)+CM_\Gamma^2\delta^{1/6},\\
 M_\Gamma^2=sup_k\|Q_k-Q\|_{W^{1,3/2}(C'_{h_0})}^2
              +\|p-c_\Gamma\|_{W^{1,3/2}(C'_{h_0})}^2.
\end{gathered}
\end{equation}
Taking first $\limsup_k$ and then $\delta\downarrow0$ proves the claim. Equivalently, for every $\varepsilon>0$ one can choose
\[
 CM_\Gamma^2\delta^{1/6}\le\varepsilon^2/2,\qquad
 C(e_k^2+n_k^2)\le\varepsilon^2/2\quad(k\ge k_0(\delta)),
\]
which yields $\|q_k-c_\Gamma\|_{L^2(\Gamma_h)}\le\varepsilon$ simultaneously for $k\ge k_0$ and almost every $0<h\le\delta$. Finite boundary-component extraction preserves all conclusions.
\end{proof}

\section{Transverse collars and boundary integrals}\label{sec:collars}

\subsection{Collar coordinates, measures, and admissible components}

\begin{lemma}[A transverse Lipschitz collar]\label{lemma:4.1}
Let $\Gamma$ be a union of connected components of $\partial\Omega$.
There exist a smooth vector field $\mathcal X$, defined near $\Gamma$,
and $h_0>0$ such that its flow
\begin{equation}\label{eq:4.1}
 \partial_h T_h(s)=\mathcal X(T_h(s)),\qquad
 T_0(s)=s,\qquad X(s,h):=T_h(s)
\end{equation}
is a bi-Lipschitz map from $\Gamma\times(0,h_0)$ onto an interior
collar $C_{h_0}$. For $0<h<h_0$, put
\[
 \Gamma_h=T_h(\Gamma),\qquad C_h=X(\Gamma\times(0,h)).
\]
There are constants $c_d,C_d,c_J,C_J>0$, independent of $h$, such that
\begin{equation}\label{eq:4.2}
 c_dh\leq\operatorname{dist}(T_h(s),\partial\Omega)\leq C_dh,
 \qquad |T_h(s)-s|\leq C_dh,
\end{equation}
\begin{equation}\label{eq:4.3}
 c_J\leq J_S(s,h)\leq C_J,\qquad
 c_J\leq J_V(s,h)\leq C_J\quad\text{a.e.},
\end{equation}
where $J_S$ and $J_V$ are the surface and volume Jacobians.
The surfaces $\Gamma_h$ have uniformly bounded Lipschitz chart
constants and satisfy
\begin{equation}\label{eq:4.4}
 \mathcal H^{n-1}(\Gamma_h\cap B(y,\rho))\leq C_A\rho^{n-1}
 \qquad (y\in\mathbb R^n,\ 0<\rho<\rho_0).
\end{equation}
If $n=3$, $\Omega$ is axisymmetric, and $\Gamma$ is separated from
the axis, the field and flow can be chosen equivariant under rotations.
The resulting collar is also separated from the axis.
\end{lemma}

\begin{proof}
Write a boundary chart, after a rigid change of coordinates, as
\[
 \Omega\cap U_\alpha
   =\{(x',x_n)\in U_\alpha:x_n>\varphi_\alpha(x')\},
 \qquad \varphi_\alpha\in W^{1,\infty}.
\]
The constant inward graph direction $e_\alpha$ satisfies
\[
 e_\alpha\cdot n_\Omega(s)
 \leq-(1+\|\nabla\varphi_\alpha\|_\infty^2)^{-1/2}
 =:-c_\alpha<0
\]
for almost every boundary point in the chart. Choose a finite subcover
and a smooth partition of unity $\chi_\alpha$ on a neighborhood of
$\Gamma$, with supports contained in these charts. Then
\begin{equation}\label{eq:4.5}
 \mathcal X=\sum_\alpha\chi_\alpha e_\alpha,\qquad
 \mathcal X(s)\cdot n_\Omega(s)\leq-c_T<0
 \quad\text{a.e. on }\Gamma.
\end{equation}

We verify the collar property for the flow of this field. For related
transverse-field and straight-displacement constructions, see
\cite{HMT,GQ}. In a smaller graph chart put
$d(x)=x_n-\varphi_\alpha(x')$. At a point of differentiability of the
graph, $\nabla d/|\nabla d|=-n_\Omega$ on the boundary.
The gradient depends only on $x'$, so smoothness of $\mathcal X$ and
\eqref{eq:4.5} give, throughout a sufficiently thin two-sided strip
in the chart,
\[
 \mathcal X(x)\cdot\nabla d(x)\geq c_T/2
 \quad\text{a.e.}
\]
For a local mollification $d_\epsilon=d*\eta_\epsilon$, the
commutator between $\mathcal X\cdot\nabla d_\epsilon$ and
$(\mathcal X\cdot\nabla d)*\eta_\epsilon$ has absolute value at
most $C\epsilon\|\nabla\mathcal X\|_\infty\|\nabla d\|_\infty$.
Integration along the smooth trajectories, followed by uniform
convergence of $d_\epsilon$ to $d$, therefore gives
\[
 d(T_t x)-d(T_u x)\geq(c_T/2)(t-u)\qquad(u<t)
\]
whenever the intervening trajectory stays in the smaller strip.
Each point sufficiently close to the graph has a unique short signed
hitting time $\tau(x)$ with $d(T_{-\tau(x)}x)=0$.
For two nearby points, strict increase and the uniform Lipschitz bound
for the short-time flow give
\[
 (c_T/2)|\tau(x)-\tau(y)|
 \leq|d(T_{-\tau(y)}x)-d(T_{-\tau(y)}y)|
 \leq C|x-y|.
\]
Hence $x\mapsto(T_{-\tau(x)}x,\tau(x))$ is the local Lipschitz
inverse of $(s,h)\mapsto T_h(s)$.
A finite cover by smaller charts and the bound
$|T_h(s)-s|\leq\|\mathcal X\|_\infty h$ make the time uniform.
If two images coincide for these times, their starting points lie in
a common chart and local uniqueness applies. Separated pairs have
separated images upon further reducing $h_0$.
This proves global bi-Lipschitz bounds onto the collar.
Since $d(T_h(s))\geq c_Th/2$ and $d$ is Lipschitz, distance to the
boundary is bounded below by a constant times $h$; displacement gives
the upper bound in \eqref{eq:4.2}.

The ambient flow satisfies
\[
 DT_h(s)=I+O(h),\qquad
 \det DT_h(s)=\exp\left(\int_0^h
       \operatorname{div}\mathcal X(T_t(s))\,dt\right).
\]
The surface area formula gives
\[
 J_S(s,h)=|\operatorname{cof}(DT_h(s))n_\Omega(s)|.
\]
The identity $DT_h(s)\mathcal X(s)=\mathcal X(T_h(s))$ gives
\[
 J_V(s,h)=|\det DT_h(s)|\,
            |\mathcal X(s)\cdot n_\Omega(s)|.
\]
Equation \eqref{eq:4.5} and the boundedness of $DT_h$ and its
inverse imply \eqref{eq:4.3}. Uniform Lipschitz charts and the area
formula imply \eqref{eq:4.4}.

For the axisymmetric assertion, let $S$ be reflection in a fixed plane
through the axis and put
\[
 \begin{split}
 \mathcal X_{\mathrm{rot}}(x)
   &=\frac1{2\pi}\int_0^{2\pi}R_{-\theta}
                    \mathcal X(R_\theta x)\,d\theta,\\
 \mathcal X_{\mathrm{ax}}(x)
   &=\frac12\bigl(\mathcal X_{\mathrm{rot}}(x)
             +S^{-1}\mathcal X_{\mathrm{rot}}(Sx)\bigr).
 \end{split}
\]
The field is first restricted to a rotation-invariant neighborhood of
$\Gamma$, so all terms in the averages are defined.
Using $n_\Omega(R_\theta s)=R_\theta n_\Omega(s)$ a.e., we have
\[
 \mathcal X_{\mathrm{rot}}(s)\cdot n_\Omega(s)
 =\frac1{2\pi}\int_0^{2\pi}
       \mathcal X(R_\theta s)\cdot n_\Omega(R_\theta s)\,d\theta
 \leq-c_T.
\]
The same inequality holds after reflection averaging, since an
axisymmetric domain is also invariant under $S$.
The averaged field is invariant under $O(2)$ and hence has zero
azimuthal component. This averaging concerns the geometric field
only; it imposes no reflection symmetry on the Navier--Stokes data.
Uniqueness for the flow equation implies
$T_h(R_\theta s)=R_\theta T_h(s)$, and the flow preserves each
meridional half-plane away from the axis.
If $d_*:=\operatorname{dist}(\Gamma,\{x_1=x_2=0\})>0$, reduce
$h_0$ so that $C_dh_0<d_*/2$; then $r\geq d_*/2$ on the collar.
\end{proof}

For every nonnegative measurable $F$, the change-of-variables formulas
are
\begin{equation}\label{eq:4.6}
 \int_{\Gamma_h}F(y)\,dS_y
   =\int_\Gamma F(T_h(s))J_S(s,h)\,dS_s,
\end{equation}
\begin{equation}\label{eq:4.7}
 \int_{C_\delta}F(x)\,dx
   =\int_0^\delta\int_\Gamma
       F(T_h(s))J_V(s,h)\,dS_s\,dh.
\end{equation}
Thus fixed-coordinate and moving-surface integrals are uniformly
equivalent. In particular,
\begin{equation}\label{eq:4.8}
 |C_\delta|\leq C\delta,\qquad
 \mathcal H^{n-1}(\Gamma_h)\leq C.
\end{equation}
All statements involving Sobolev restrictions initially hold for
almost every $h$. Subsequent choices of $h$ are made in the
intersection of the full-measure sets required by these statements.

A collar is called admissible if $n=2$, or if $n=3$ and its boundary
components are separated from the axis. In the latter case its
meridional coordinates satisfy $0<r_*\leq r\leq r^*<\infty$, and
for axisymmetric scalar integrands,
\begin{equation}\label{eq:4.9}
 \begin{split}
 \int_{C_\delta}F(r,z)\,dx
   &=2\pi\int_{\widehat C_\delta}F(r,z)r\,dr\,dz,\\
 \int_{\Gamma_h}F(r,z)\,dS
   &=2\pi\int_{\gamma_h}F(r,z)r\,d\mathcal H^1.
 \end{split}
\end{equation}
Here $\widehat C_\delta$ and $\gamma_h$ denote the meridional
collar and its generating curve. Weighted and unweighted norms in
these formulas are equivalent, uniformly in $h$.

\subsection{Uniform Sobolev traces and a harmonic pressure estimate}

\begin{lemma}[Sobolev and harmonic collar estimates]\label{lemma:4.2}
For every $b\in H^1(\Omega)$,
\begin{equation}\label{eq:4.10}
 \operatorname*{ess\,sup}_{0<h<h_0}
       \|b\|_{L^2(\Gamma_h)}\leq C\|b\|_{H^1(\Omega)}.
\end{equation}
If $\pi\in L^2(\Omega)$ and $\Delta\pi=0$ in distributions, then
\begin{equation}\label{eq:4.11}
 \|\pi\|_{L^2(\Gamma_h)}^2
       \leq Ch^{-1}\|\pi\|_{L^2(\Omega)}^2
       \qquad(0<h<h_0).
\end{equation}
For the remainder $r_k=(\mu_k/J_k)\pi$ in \eqref{eq:2.12},
\begin{equation}\label{eq:4.12}
 \sup_{\mu_k^2\leq h<h_0}\|r_k\|_{L^2(\Gamma_h)}
       \leq C\mu_k\longrightarrow0.
\end{equation}
\end{lemma}

\begin{proof}
For smooth $b$, put $B(s,h)=b(T_h(s))$ and
$M_{\mathcal X}=\sup|\mathcal X|$. Then
\[
 \partial_hB(s,h)=\nabla b(T_h(s))\mathcal X(T_h(s)),\qquad
 |\partial_h|B(s,h)|^2|
 \leq2M_{\mathcal X}|b(T_h(s))|\,|\nabla b(T_h(s))|.
\]
Integration first in $h$ and then on $\Gamma$ yields
\[
 \begin{split}
 \int_\Gamma |B(s,h)|^2\,dS_s
 &\leq\|\operatorname{Tr}b\|_{L^2(\Gamma)}^2
       +2M_{\mathcal X}\int_0^h\int_\Gamma
          |b(T_t(s))|\,|\nabla b(T_t(s))|\,dS_s\,dt\\
 &\leq C\|b\|_{H^1(\Omega)}^2
       +C\int_{C_{h_0}}|b|\,|\nabla b|\,dx\\
 &\leq C\|b\|_{H^1(\Omega)}^2.
 \end{split}
\]
Use \eqref{eq:4.6} to obtain \eqref{eq:4.10}.
Approximation by smooth ambient functions, available by the $H^1$
extension theorem on a Lipschitz domain, gives the assertion for
$b\in H^1$. Equivalently, pullback to the collar gives
$B\in H^1((0,h_0);L^2(\Gamma))$, and the preceding calculation
uses its continuous $L^2(\Gamma)$ representative.

For harmonic $\pi$, interior regularity permits pointwise evaluation
in $\Omega$. Fix $\alpha\in(0,c_d/2)$, so that
$B(y,\alpha h)\subset\Omega$ for $y\in\Gamma_h$.
Subharmonicity of $|\pi|^2$ gives
\[
 |\pi(y)|^2\leq Ch^{-n}\int_{B(y,\alpha h)}|\pi(x)|^2\,dx.
\]
Integrating and applying Tonelli's theorem,
\[
 \begin{split}
 \int_{\Gamma_h}|\pi(y)|^2\,dS_y
 &\leq Ch^{-n}\int_\Omega|\pi(x)|^2
       \left(\int_{\Gamma_h}
          \mathbf 1_{\{|x-y|<\alpha h\}}\,dS_y\right)dx\\
 &\leq Ch^{-n}\int_\Omega|\pi(x)|^2
       \mathcal H^{n-1}(\Gamma_h\cap B(x,\alpha h))\,dx\\
 &\leq Ch^{-1}\|\pi\|_{L^2(\Omega)}^2.
 \end{split}
\]
The final inequality follows from \eqref{eq:4.4}, after reducing
$h_0$ if necessary so that $\alpha h_0<\rho_0$.
Since $J_k^{-1}=\lambda_k\mu_k/\nu$ and
$0\leq\lambda_k\leq1$,
\[
 \begin{split}
 \|r_k\|_{L^2(\Gamma_h)}
 &\leq C\frac{\mu_k}{J_k}h^{-1/2}\|\pi\|_2
   =C\frac{\lambda_k}{\nu}\mu_k^2h^{-1/2}\|\pi\|_2\\
 &\leq C\mu_k\qquad(h\geq\mu_k^2).
 \end{split}
\]
This proves \eqref{eq:4.12}.
\end{proof}

\subsection{Smallness of the actual head on a truncated collar}

Write
\[
 K_k:=\tfrac12|u_k|^2,\qquad
 \Phi_k=p_k+K_k=q_k+r_k+K_k.
\]
On an admissible two-dimensional collar, including a meridional
collar, the energy bound gives
\begin{equation}\label{eq:4.13}
 \|K_k\|_{W^{1,3/2}}
 \leq C\bigl(\|u_k\|_{L^3}^2
                  +\|u_k\|_{L^6}\|\nabla u_k\|_{L^2}\bigr)
 \leq C,
\end{equation}
\begin{equation}\label{eq:4.14}
 \|K_k\|_{L^6}\leq C\|K_k\|_{W^{1,3/2}}\leq C,
 \qquad\|K_k\nabla K_k\|_{L^{6/5}}\leq C.
\end{equation}
In cylindrical coordinates, derivatives of the vector components have
bounded coefficients on an admissible collar; hence the same
estimates apply to $|u_k|^2$ there.

The trace product identity follows by $H^1$ approximation and the
boundary embedding into $L^4$:
\begin{equation}\label{eq:4.15}
 \operatorname{Tr}K_k=\frac1{2J_k^2}|a|^2,\qquad
 \|\operatorname{Tr}K_k\|_{L^2(\Gamma)}^2
       =\frac1{4J_k^4}\|a\|_{L^4(\Gamma)}^4.
\end{equation}
For $n=3$, $H^{1/2}(\partial\Omega)\hookrightarrow L^4(\partial\Omega)$;
for $n=2$, the same embedding follows from the one-dimensional
fractional Sobolev theorem. The constants in these embeddings depend
on the fixed boundary.

Apply the flow calculation to $K_k^2$, now using \eqref{eq:4.14}
and H\"older's inequality with exponents $6/5$ and $6$:
\begin{equation}\label{eq:4.16}
 \begin{split}
 \operatorname*{ess\,sup}_{0<h\leq\delta}
      \|K_k\|_{L^2(\Gamma_h)}^2
 &\leq C\|\operatorname{Tr}K_k\|_{L^2(\Gamma)}^2
      +C\int_{\widehat C_\delta}|K_k|\,|\nabla K_k|\,dx'\\
 &\leq CJ_k^{-4}\|a\|_{L^4(\Gamma)}^4
      +C|\widehat C_\delta|^{1/6}
                    \|K_k\nabla K_k\|_{L^{6/5}}\\
 &\leq CJ_k^{-4}+C\delta^{1/6}.
 \end{split}
\end{equation}
In the planar case, the hat denotes the actual collar.
In the spatial case, \eqref{eq:4.9} transfers this estimate to the
physical surface.

For an admissible component $\Gamma_j$, Lemma~\ref{lemma:3.5} gives
\begin{equation}\label{eq:4.17}
 \lim_{\delta\downarrow0}\limsup_{k\to\infty}
 \operatorname*{ess\,sup}_{0<h\leq\delta}
       \|q_k-p_j\|_{L^2(\Gamma_{j,h})}=0.
\end{equation}
Combining \eqref{eq:4.12}, \eqref{eq:4.16}, and \eqref{eq:4.17},
\begin{equation}\label{eq:4.18}
 \begin{split}
 \operatorname*{ess\,sup}_{\mu_k^2\leq h\leq\delta}
       \|\Phi_k-p_j\|_{L^2(\Gamma_{j,h})}
 &\leq\operatorname*{ess\,sup}_{0<h\leq\delta}
       \|q_k-p_j\|_{L^2(\Gamma_{j,h})}\\
 &\hspace{1em}+C\mu_k+CJ_k^{-2}+C\delta^{1/12}.
 \end{split}
\end{equation}
Thus, for a finite union $\Gamma$ of admissible components with
$p_j=0$, the precise quantifiers are
\begin{equation}\label{eq:4.19}
 \begin{gathered}
 \forall\sigma>0\quad\exists\delta_\sigma\in(0,h_0),\ k_\sigma:
 \\
 \|\Phi_k\|_{L^2(\Gamma_h)}\leq\sigma
 \quad\text{for every }k\geq k_\sigma
 \text{ and a.e. }h\in[\mu_k^2,\delta_\sigma].
 \end{gathered}
\end{equation}
Indeed, first choose $\delta_\sigma$ using the outer limit in
\eqref{eq:4.17} and the term $C\delta^{1/12}$; then choose
$k_\sigma$ using the inner limit and
$\mu_k+J_k^{-2}\to0$. Increase $k_\sigma$ so that
$\mu_k^2<\delta_\sigma$.

\subsection{Simultaneous selection of the two velocity bounds}

\begin{lemma}[A common surface for the velocity estimates]\label{lemma:4.3}
For every $\epsilon>0$ there is $C_\epsilon<\infty$, independent
of $k$ and of $\delta$, such that for each fixed $0<\delta<h_0$
and all sufficiently large $k$, there exists
$h_k\in[\mu_k^2,\delta]$ satisfying
\begin{equation}\label{eq:4.20}
 \int_{\Gamma_{h_k}}|u_k|^2\,dS\leq C_\epsilon\mu_k^2,
 \qquad
 \int_{\Gamma_{h_k}}|u_k|\,|\nabla u_k|\,dS<\epsilon.
\end{equation}
The surface may be selected in any prescribed full-measure set of
collar parameters. In particular, it may satisfy \eqref{eq:4.19}
simultaneously with \eqref{eq:4.20}.
\end{lemma}

\begin{proof}
Fix $k$ and abbreviate $\mu=\mu_k$. Integrals in the next
calculation use the fixed measure on $\Gamma$:
\[
 \begin{aligned}
 g(h)&=\int_\Gamma|u_k(T_h(s))|^2\,dS_s,
 &G(h)&=g(h)+\mu^2,\\
 F(h)&=\int_\Gamma|u_k(T_h(s))|\,
                    |\nabla u_k(T_h(s))|\,dS_s,
 &E(h)&=\int_\Gamma|\nabla u_k(T_h(s))|^2\,dS_s.
 \end{aligned}
\]
The energy bound and \eqref{eq:4.7} give
\begin{equation}\label{eq:4.21}
 \int_0^{h_0}E(h)\,dh\leq C_*,\qquad
 |G'(h)|\leq C_{\mathcal X}F(h),\qquad
 F(h)^2\leq g(h)E(h)\leq G(h)E(h),
\end{equation}
where $C_{\mathcal X}=2M_{\mathcal X}$ may be increased to a fixed
positive constant. The function $G$ is absolutely continuous and
$G\geq\mu^2$.

The trace and the fundamental theorem along the flow give
\[
 \begin{split}
 g(0)&=J_k^{-2}\|a\|_{L^2(\Gamma)}^2\leq C_a\mu^2,\\
 |u_k(T_h(s))|^2
 &\leq2|J_k^{-1}a(s)|^2
      +2M_{\mathcal X}^2h
                \int_0^h|\nabla u_k(T_t(s))|^2\,dt,\\
 g(h)&\leq2g(0)+2M_{\mathcal X}^2h\int_0^hE(t)\,dt.
 \end{split}
\]
Consequently, with $h_a=\mu^2$,
\begin{equation}\label{eq:4.22}
 G(h_a)\leq C_0\mu^2,\qquad
 C_0:=1+2C_a+2M_{\mathcal X}^2C_*.
\end{equation}
Let $C_J$ be the upper surface-Jacobian bound in \eqref{eq:4.3}.
Choose
\begin{equation}\label{eq:4.23}
 \eta:=\frac\epsilon{2C_J},\qquad
 L:=C_0\exp\left(1+\frac{C_{\mathcal X}C_*}{\eta}\right),
 \qquad C_\epsilon:=C_JL.
\end{equation}
These constants do not depend on $\delta$ or $k$.

\emph{Case 1: a first crossing occurs.}
Suppose $G(h_1)=L\mu^2$ for a first $h_1\in(h_a,\delta]$.
Then $G\leq L\mu^2$ on $[h_a,h_1]$.
If $F\geq\eta$ almost everywhere on this interval,
\[
 \begin{split}
 \log\frac L{C_0}
 &\leq\log\frac{G(h_1)}{G(h_a)}
    =\int_{h_a}^{h_1}\frac{G'(h)}{G(h)}\,dh\\
 &\leq C_{\mathcal X}\int_{h_a}^{h_1}\frac{F(h)}{G(h)}\,dh
  \leq\frac{C_{\mathcal X}}\eta
           \int_{h_a}^{h_1}\frac{F(h)^2}{G(h)}\,dh\\
 &\leq\frac{C_{\mathcal X}}\eta
           \int_{h_a}^{h_1}E(h)\,dh
  \leq\frac{C_{\mathcal X}C_*}\eta.
 \end{split}
\]
This contradicts \eqref{eq:4.23}. Hence
$\{h\in(h_a,h_1):F(h)<\eta\}$ has positive measure, and
$g(h)\leq L\mu^2$ throughout that set.

\emph{Case 2: no crossing occurs.}
Then $G(h)<L\mu^2$ for $h\in[h_a,\delta]$, and
\[
 \begin{split}
 \int_{h_a}^{\delta}F(h)\,dh
 &\leq\left(\int_{h_a}^{\delta}G(h)\,dh\right)^{1/2}
       \left(\int_{h_a}^{\delta}E(h)\,dh\right)^{1/2}\\
 &\leq\sqrt{LC_*}\,\mu\sqrt{\delta-h_a}.
 \end{split}
\]
For fixed $\delta$, choose $k$ sufficiently large that
\[
 \mu^2<\frac\delta2,\qquad
 \sqrt{\frac{2LC_*}\delta}\,\mu<\frac\eta2.
\]
The average of $F$ on $[h_a,\delta]$ is then less than
$\eta/2$. Therefore $\{F<\eta\}$ has positive measure on that
interval, where again $g\leq L\mu^2$.

In both cases choose $h_k$ in this positive-measure set and in the
intersection of the required full-measure sets.
By \eqref{eq:4.6},
\[
 \begin{split}
 \int_{\Gamma_{h_k}}|u_k|^2\,dS
   &\leq C_Jg(h_k)\leq C_\epsilon\mu_k^2,\\
 \int_{\Gamma_{h_k}}|u_k|\,|\nabla u_k|\,dS
   &\leq C_JF(h_k)<\epsilon/2<\epsilon.
 \end{split}
\]
This proves all assertions.
\end{proof}

\subsection{Interior regularity for the individual normalized solutions}

All estimates in this subsection are for a fixed $k$; their constants
may depend on $\mu_k^{-1}$.
Write $u=u_k$, $p=p_k$, $g=g_k$, and $\mu=\mu_k>0$.
Interior Stokes estimates \cite{Galdi}, on successively smaller
relatively compact open sets, give the following bootstrap in
dimension three:
\begin{equation}\label{eq:4.24}
 \begin{aligned}
 u\in H^1,\quad g\in H^1
   &\ \Longrightarrow\ g-u\cdot\nabla u\in L^{3/2}\\
   &\ \Longrightarrow\ u\in W^{2,3/2}_{\mathrm{loc}},
                      \quad p\in W^{1,3/2}_{\mathrm{loc}},\\
 \nabla u\in L^3_{\mathrm{loc}},\quad u\in L^6
   &\ \Longrightarrow\ u\cdot\nabla u\in L^2_{\mathrm{loc}}\\
   &\ \Longrightarrow\ u\in H^2_{\mathrm{loc}},
                      \quad p\in H^1_{\mathrm{loc}}.
 \end{aligned}
\end{equation}
Now $u\in L^\infty_{\mathrm{loc}}$ and
$\nabla u\in L^6_{\mathrm{loc}}$.
On a compact interior set $K$,
\[
 \|\nabla(u\cdot\nabla u)\|_{L^2(K)}
 \leq C\|\nabla u\|_{L^4(K)}^2
       +C\|u\|_{L^\infty(K)}\|D^2u\|_{L^2(K)}<\infty.
\]
A further interior Stokes estimate therefore gives
\begin{equation}\label{eq:4.25}
 u_k\in H^3_{\mathrm{loc}}(\Omega),\qquad
 p_k\in H^2_{\mathrm{loc}}(\Omega),\qquad
 \Phi_k\in H^2_{\mathrm{loc}}(\Omega).
\end{equation}
The planar bootstrap follows from the same estimates and the
stronger two-dimensional Sobolev embeddings.
In the spatial axisymmetric case, restriction to a compact
meridional set with $r>0$ gives the corresponding two-dimensional
Sobolev regularity. In particular, the Sobolev Morse--Sard theorem
is applied to the planar or meridional $H^2$ head, not to an
arbitrary scalar function on a three-dimensional domain.

\subsection{The head equation, with all force terms retained}

In the plane let
\[
 \omega=\partial_1u_2-\partial_2u_1,\qquad
 v^\perp=(-v_2,v_1),\qquad
 \nabla^\perp\omega=(-\partial_2\omega,\partial_1\omega).
\]
Using $\operatorname{div}u=0$, direct differentiation yields
\[
 \begin{aligned}
 \Delta u&=\nabla^\perp\omega,\\
 (u\cdot\nabla)u&=\nabla(|u|^2/2)+\omega u^\perp,\\
 \nabla\Phi&=\mu\nabla^\perp\omega-\omega u^\perp+g.
 \end{aligned}
\]
Since $\operatorname{div}u^\perp=-\omega$ and
$\nabla\omega\cdot u^\perp=-u\cdot\nabla^\perp\omega$,
\[
 \begin{aligned}
 \Delta\Phi
   &=-\operatorname{div}(\omega u^\perp)+\operatorname{div}g
     =\omega^2+u\cdot\nabla^\perp\omega+\operatorname{div}g,\\
 u\cdot\nabla\Phi&=\mu u\cdot\nabla^\perp\omega+g\cdot u.
 \end{aligned}
\]
Eliminating $u\cdot\nabla^\perp\omega$ gives the desired head
equation in the plane.

In three dimensions let $\omega=\operatorname{curl}u$.
The vector identities
\[
 \begin{aligned}
 \Delta u&=-\operatorname{curl}\omega,\\
 (u\cdot\nabla)u&=\nabla(|u|^2/2)-u\times\omega,\\
 \nabla\Phi&=u\times\omega-\mu\operatorname{curl}\omega+g,\\
 \operatorname{div}(u\times\omega)
      &=|\omega|^2-u\cdot\operatorname{curl}\omega
 \end{aligned}
\]
give
\[
 \begin{aligned}
 \Delta\Phi
   &=|\omega|^2-u\cdot\operatorname{curl}\omega
                         +\operatorname{div}g,\\
 u\cdot\nabla\Phi
   &=-\mu u\cdot\operatorname{curl}\omega+g\cdot u.
 \end{aligned}
\]
Thus, in either dimension,
\begin{equation}\label{eq:4.26}
 \Delta\Phi_k=|\omega_k|^2
      +\mu_k^{-1}\operatorname{div}(\Phi_ku_k)
      -\mu_k^{-1}g_k\cdot u_k+\operatorname{div}g_k.
\end{equation}
Here $u_k\cdot\nabla\Phi_k=\operatorname{div}(\Phi_ku_k)$
follows from $\operatorname{div}u_k=0$.
No condition $\operatorname{div}f=0$ is imposed on the regular
force considered in this section.

\subsection{Normal head flux and uniform force errors}

Let $n_h$ be either consistent orientation of a connected closed
collar surface. In dimension two, if $\tau_h$ is the corresponding
unit tangent, then
\[
 \int_{\Gamma_h}\nabla^\perp\omega_k\cdot n_h\,dS
    =\pm\int_{\Gamma_h}\partial_{\tau_h}\omega_k\,dS=0.
\]
In dimension three, the intrinsic Stokes formula on a closed
oriented Lipschitz surface gives
\[
 \int_{\Gamma_h}\operatorname{curl}\omega_k\cdot n_h\,dS=0.
\]
To see why an enclosed hole causes no difficulty, restrict the
one-form $\omega_k\cdot dx$ to the surface.
The integrand is the surface pullback of its exterior derivative,
and
\[
 \int_{\Gamma_h}d(\omega_k\cdot dx)
     =\int_{\partial\Gamma_h}\omega_k\cdot dx=0.
\]
This calculation requires $\omega_k$ only in a neighborhood of
$\Gamma_h$. Smooth approximation in that neighborhood justifies
the identity at the regularity \eqref{eq:4.25}.
It does not extend the fluid field through the hole.

The first-order head formulas therefore imply, in both dimensions
and also for a finite union of components,
\begin{equation}\label{eq:4.27}
 \begin{split}
 \left|\int_{\Gamma_h}\nabla\Phi_k\cdot n_h\,dS\right|
 &\leq C\int_{\Gamma_h}|u_k|\,|\nabla u_k|\,dS
      +\left|\int_{\Gamma_h}g_k\cdot n_h\,dS\right|\\
 &\leq C\int_{\Gamma_h}|u_k|\,|\nabla u_k|\,dS
      +\mathcal H^{n-1}(\Gamma_h)^{1/2}
                             \|g_k\|_{L^2(\Gamma_h)}\\
 &\leq C\int_{\Gamma_h}|u_k|\,|\nabla u_k|\,dS+C\mu_k^2.
 \end{split}
\end{equation}
The last line uses \eqref{eq:4.10}, \eqref{eq:4.8}, and
$\|g_k\|_{H^1}=O(\mu_k^2)$.
On the selected surface,
\begin{equation}\label{eq:4.28}
 \left|\int_{\Gamma_{h_k}}\nabla\Phi_k\cdot n_h\,dS\right|
       \leq C\epsilon+C\mu_k^2.
\end{equation}

For every measurable set $D\subset\Omega$, the force terms satisfy
bounds independent of the particular set:
\begin{equation}\label{eq:4.29}
 \begin{split}
 \mu_k^{-1}\left|\int_D g_k\cdot u_k\,dx\right|
   &\leq\mu_k^{-1}\|g_k\|_{L^2(\Omega)}
                       \|u_k\|_{L^2(\Omega)}\leq C\mu_k,\\
 \left|\int_D\operatorname{div}g_k\,dx\right|
   &\leq|\Omega|^{1/2}
                    \|\operatorname{div}g_k\|_{L^2(\Omega)}
     \leq C\mu_k^2.
 \end{split}
\end{equation}
Finally, whenever $D\Subset\Omega$ and $D$ is bounded by a
collar surface and finitely many interior regular level surfaces,
\eqref{eq:4.25} and Gauss--Green give the explicit integrated identity
\begin{equation}\label{eq:4.30}
 \begin{split}
 \int_{\partial D}\nabla\Phi_k\cdot n_D\,dS
   &=\int_D|\omega_k|^2\,dx
      +\mu_k^{-1}\int_{\partial D}\Phi_ku_k\cdot n_D\,dS\\
   &\hspace{1em}-\mu_k^{-1}\int_Dg_k\cdot u_k\,dx
      +\int_D\operatorname{div}g_k\,dx.
 \end{split}
\end{equation}
Together, \eqref{eq:4.19}, \eqref{eq:4.20}, \eqref{eq:4.27}, and
\eqref{eq:4.29} are the boundary inputs for the separating-level
argument.

\section{Euler geometry and the contradiction}\label{sec:euler}

Throughout this section the normalized fields satisfy
\begin{equation}\label{eq:5.1}
\begin{gathered}
 u_k\rightharpoonup v\quad\hbox{in }H^1(\Omega),\qquad
 u_k\longrightarrow v\quad\hbox{in }L^4(\Omega),\qquad
 \operatorname{Tr}v=0,\\
 p_k\longrightarrow p\quad\hbox{in }L^2(\Omega),\qquad
 \nabla p=-(v\cdot\nabla)v,\qquad \operatorname{div}v=0,\\
 \mu_k\downarrow0,\qquad \lambda_k\longrightarrow\lambda_0>0,\qquad
 \|g_k\|_{H^1(\Omega)}\le C_g\mu_k^2.
\end{gathered}
\end{equation}
The last estimate is used when integrating the viscous head equation for the
regular-force sequence. The limiting Euler identities, convergence on regular
cycles, and separating sets use only the velocity, pressure, and interior head
estimates. Section~\ref{sec:dual} verifies these inputs for its corrected
sequence and estimates the effective force separately. Extracted subsequences
are relabelled. Set
\begin{equation}\label{eq:5.2}
 \Phi=p+\frac{|v|^2}{2},\qquad
 \Phi_k=p_k+\frac{|u_k|^2}{2},\qquad
 \omega=
 \begin{cases}
 \partial_1v_2-\partial_2v_1,&n=2,\\
 \operatorname{curl}v,&n=3.
 \end{cases}
\end{equation}
In particular,
\begin{equation}\label{eq:5.3}
\begin{gathered}
 \big\||u_k|^2-|v|^2\big\|_2
 \le(\|u_k\|_4+\|v\|_4)\|u_k-v\|_4\longrightarrow0,\\
 \Phi_k\longrightarrow\Phi\quad\hbox{in }L^2(\Omega),\qquad
 \sup_k\|\Phi_k\|_{W^{1,3/2}(K)}<\infty
 \quad(K\Subset\Omega).
\end{gathered}
\end{equation}
For spatial fields the last estimate also holds in ordinary two-dimensional
coordinates on every compact meridional set with \(r>0\). The Euler identities are
\begin{equation}\label{eq:5.4}
 \nabla\Phi=
 \begin{cases}
 -\omega v^\perp,&n=2,\quad v^\perp=(-v_2,v_1),\\
 v\times\omega,&n=3.
 \end{cases}
\end{equation}

\subsection{Boundary constants and the energy identity}

\begin{lemma}[Euler boundary pressures]\label{lemma:5.1}
There are \(p_0,\ldots,p_N\in\mathbb R\) such that
\begin{equation}\label{eq:5.5}
 \operatorname{Tr}_{\Gamma_j}p
 =\operatorname{Tr}_{\Gamma_j}\Phi=p_j
 \quad\mathcal H^{n-1}\hbox{-a.e. on }\Gamma_j.
\end{equation}
In the axisymmetric case, all components meeting the axis have the same
constant \(p_{\rm ax}\). Moreover,
\begin{equation}\label{eq:5.6}
 \nu=-\lambda_0\sum_{j=0}^N p_jF_j,\qquad
 F_j=\int_{\Gamma_j}a\cdot n_\Omega\,dS.
\end{equation}
Consequently, \(\max_jp_j>\min_jp_j\).
\end{lemma}

\begin{proof}
We first verify constancy directly on the physical Lipschitz boundary; compare
\cite[Section~3]{KPR2013} and \cite[Theorem~3.1]{KPR2015axis}. Let
\(d(x)=\operatorname{dist}(x,\partial\Omega)\). The zero-trace Hardy inequality
and the Euler equation give
\begin{equation}\label{eq:5.7}
 \int_\Omega\frac{|\nabla p|}{d}\,dx
 \le\|v/d\|_2\|\nabla v\|_2
 \le C_\Omega\|\nabla v\|_2^2.
\end{equation}
For a boundary component use Lemma~\ref{lemma:4.1}; equivariance is not needed
here. Its Jacobian and distance bounds imply
\[
 \int_0^{h_0}\frac1h\int_\Gamma
       |\nabla p(T_hs)|\,dS_s\,dh<\infty.
\]
There is consequently a sequence of admissible slices \(h_\ell\downarrow0\)
on which the inner integral tends to zero. In boundary charts the tangential
weak derivatives of \(p\circ T_{h_\ell}\) are bounded in \(L^1\) by a fixed
multiple of this integral. Also
\(p\circ T_h\to\operatorname{Tr}p\) in \(L^1(\Gamma)\): pull back
\(p\in W^{1,3/2}\) to the collar and integrate its \(h\)-derivative. Passing
to weak tangential derivatives against smooth chart tests shows that
\(\operatorname{Tr}p\) is constant on each chart. Connectedness identifies
these constants on \(\Gamma\).

Equality on the axis-intersecting components is the axisymmetric Euler
pressure theorem \cite[Theorem~3.3]{KPR2015axis}, applied with its parameter
equal to one and
\[
 v\in H_0^1(\Omega),\qquad p\in W^{1,3/2}(\Omega),\qquad
 (v\cdot\nabla)v+\nabla p=0.
\]
In the plane one also has \(p\in W^{1,q}\) for every \(1<q<2\).
The Sobolev chain rule and trace continuity give
\[
 |v|^2\in W^{1,3/2}(\Omega),\qquad
 \operatorname{Tr}(|v|^2)=|\operatorname{Tr}v|^2=0,
\]
which proves the assertion for the head.

Let \(U\) be the fixed solenoidal Stokes lift used in the continuation
argument. The products needed for Gauss--Green satisfy
\begin{equation}\label{eq:5.8}
\begin{aligned}
 \|pU\|_2&\le\|p\|_3\|U\|_6,\\
 \|U\otimes\nabla p\|_{6/5}&\le\|U\|_6\|\nabla p\|_{3/2},\\
 \|p\nabla U\|_{6/5}&\le\|p\|_3\|\nabla U\|_2.
\end{aligned}
\end{equation}
These exponents suffice in both dimensions. Hence \(pU\in W^{1,6/5}\), and
the product trace is
\[
 \operatorname{Tr}(pU)=(\operatorname{Tr}p)a\in L^{4/3}(\partial\Omega),
 \qquad
 \|\operatorname{Tr}(pU)\|_{4/3}
 \le\|\operatorname{Tr}p\|_2\|a\|_4.
\]
Since \(\operatorname{div}U=0\),
\begin{equation}\label{eq:5.9}
\begin{aligned}
 b(v,v,U)
 &=\int_\Omega(v\cdot\nabla)v\cdot U\,dx\\
 &=-\int_\Omega\nabla p\cdot U\,dx
   =-\int_\Omega\operatorname{div}(pU)\,dx\\
 &=-\int_{\partial\Omega}(\operatorname{Tr}p)a\cdot n_\Omega\,dS
   =-\sum_{j=0}^Np_jF_j.
\end{aligned}
\end{equation}
Equation~\eqref{eq:2.8} and strong \(L^4\) convergence give
\(\nu=\lambda_0b(v,v,U)\), proving \eqref{eq:5.6}. If all \(p_j=c\),
then \(\nu=-\lambda_0c\sum_jF_j=0\), a contradiction.
\end{proof}

We may shift all pressure-related quantities simultaneously:
\begin{equation}\label{eq:5.10}
 (p_k,q_k,p,\Phi_k,\Phi,p_j)
 \longmapsto
 (p_k-c,q_k-c,p-c,\Phi_k-c,\Phi-c,p_j-c).
\end{equation}
All differential identities, convergences, and \eqref{eq:5.6} are unchanged,
because \(\sum_jF_j=0\). This shift is made only for the normalized geometric
argument; it does not change the final mean-zero normalization of a pressure
associated with a fixed force representative.

\subsection{Streamfunctions, Euler continuation, and local cycle convergence}

In the plane, let \(\widetilde v\) be the zero extension. For
\(\varphi\in C_c^\infty(\mathbb R^2)\),
\(\int_{\mathbb R^2}\widetilde v\cdot\nabla\varphi
=\int_\Omega v\cdot\nabla\varphi=0\). Thus
\begin{equation}\label{eq:5.11}
 \widetilde v\in H^1(\mathbb R^2),\qquad
 \nabla\psi=\widetilde v^\perp,\qquad
 \psi\in H^2_{\rm loc}(\mathbb R^2),\qquad
 \nabla\psi=0\quad\hbox{in }\mathbb R^2\setminus\overline\Omega.
\end{equation}
The continuous representative of \(\psi\) is constant on each complementary
component and its boundary.

In space put
\begin{equation}\label{eq:5.12}
 D=\{(r,z):r>0,\ (r,0,z)\in\Omega\},\qquad
 \partial_r\psi=-rv_z,\qquad \partial_z\psi=rv_r.
\end{equation}
The compatibility relation is
\[
 \partial_z(-rv_z)=\partial_r(rv_r)
 \quad\Longleftrightarrow\quad
 \partial_r(rv_r)+\partial_z(rv_z)=0.
\]
For a connected Lipschitz subdomain \(D_0\Subset D\),
\[
 \inf_{c\in\mathbb R}\|\psi-c\|_{H^2(D_0)}
 \le C_{D_0}\|v\|_{H^1(\widetilde D_0)},
\]
where \(\widetilde D_0\) denotes rotation about the axis. To obtain a
single-valued streamfunction, extend \(v\) by zero to \(\mathbb R^3\) and
integrate \((-rv_z,rv_r)\) in the simply connected half-plane \(r>0\).
The extension is in \(H^1(\mathbb R^3)\), is solenoidal, and obeys the
compatibility relation distributionally across every physical wall.
Hence \(\psi\in H^2_{\rm loc}(\{r>0\})\). The meridional generators of
boundary components are connected
\cite[condition (S1), Section~3]{KPR2015axis}; the extended \(\psi\) is
constant on each complementary region off the axis. This is the continuation
used in \cite[Remark~3.6]{KPR2015axis}.

\begin{lemma}[Continuation of the Euler fields]\label{lemma:5.continuation}
Extend \(v\) by zero and extend \(p\) by \(p_j\) in the complementary
component adjacent to \(\Gamma_j\), including the unbounded component
when \(j=0\). Then
\[
 \widetilde p\in W^{1,3/2}_{\rm loc}(\mathbb R^n),\qquad
 \nabla\widetilde p=\mathbf1_\Omega\nabla p,
\]
and \((\widetilde v,\widetilde p)\) solves the stationary Euler equations
in \(\mathbb R^n\). One can place these fields in a smooth bounded auxiliary
domain with finitely many holes, whose boundary components lie in regions
where the extended velocity vanishes and the extended pressure is constant.
\end{lemma}

\begin{proof}
Equality of the pressure traces in \eqref{eq:5.5} and Sobolev gluing across
Lipschitz graphs give the asserted extension and derivative. Since
\(\widetilde v\in H^1\) is the zero extension and is solenoidal, the Euler
equation holds a.e.\ and distributionally after gluing.
In the plane, put a small closed disk inside each bounded hole. In space,
put a small closed ball centered on the axis inside each axis-meeting hole
and a small closed solid torus about one circular orbit in each off-axis
hole. An axis-meeting hole contains an axial interval, by the local graph
argument in Lemma~\ref{lemma:5.axis}. Enclose \(\overline\Omega\) in a
larger disk or, in space, an axis-centered ball, and remove the chosen cores.
The resulting auxiliary domain is smooth. Every core and the outer boundary
lie strictly inside regions where the fields are constant, so the extended
fields satisfy the homogeneous Euler hypotheses used in
\cite[Sections~3 and 4]{KPR2015}.
\end{proof}

The continuation in Lemma~\ref{lemma:5.continuation} is applied only to the
Euler fields. Every convergence assertion for \(\Phi_k\), and every viscous
flux identity below, is restricted to the original physical domain.

\begin{lemma}[Local cycle convergence]\label{lemma:5.2}
There is a single subsequence such that, outside a set of streamfunction
values of one-dimensional measure zero,
\begin{equation}\label{eq:5.13}
 \|\Phi_k-\Phi\|_{C(S)}\longrightarrow0
\end{equation}
on every regular interior cycle \(S\) of \(\psi\). In space the statement is
on meridional cycles separated from the axis.
\end{lemma}

\begin{proof}
Fix \(q=3/2\). First use a regular \(C^1\) foliation chart
\[
 \Xi:I\times J\longrightarrow K\Subset\Omega,\qquad
 \psi(\Xi(s,t))=t,
\]
with bounded chart derivatives and Jacobian bounded away from zero. The
same notation is used for a compact meridional chart. Set
\[
\begin{gathered}
 Z_k=\Phi_k-\Phi,\qquad z_k(s,t)=Z_k(\Xi(s,t)),\\
 A_k(t)=\|z_k(\cdot,t)\|_{L^q(I)},\quad
 B_k(t)=\|\partial_s z_k(\cdot,t)\|_{L^q(I)}.
\end{gathered}
\]
For \(z\in W^{1,q}(I)\), choose \(s_0\) such that
\(|z(s_0)|^q\le |I|^{-1}\|z\|_q^q\). The fundamental theorem for absolutely
continuous functions and H\"older's inequality give
\begin{equation}\label{eq:5.14}
 |z(s)|^q
 \le |I|^{-1}\|z\|_q^q
       +q\int_I|z|^{q-1}|z'|\,ds
 \le |I|^{-1}\|z\|_q^q+q\|z\|_q^{q-1}\|z'\|_q.
\end{equation}
Integrating in \(t\) yields
\begin{equation}\label{eq:5.15}
\begin{aligned}
 \int_J\|z_k(\cdot,t)\|_\infty^q\,dt
 &\le C\int_J A_k^q\,dt+C\int_J A_k^{q-1}B_k\,dt\\
 &\le C\|Z_k\|_{L^q(K)}^q
   +C\|Z_k\|_{L^q(K)}^{q-1}\|Z_k\|_{W^{1,q}(K)}
 \longrightarrow0.
\end{aligned}
\end{equation}
Choose \(k_\ell\) so that the left-hand side is at most \(2^{-\ell}\).
Tonelli's theorem then gives convergence to zero in \(L^\infty(I)\) for
almost every leaf \(t\).

For a Sobolev streamfunction, exhaust the interior by countably many
relatively compact patches, off the axis in a meridian. On each patch
multiply \(\psi\) by a smooth interior cutoff equal to one near its closure,
and extend by zero to the plane. The resulting function is compactly
supported and belongs to \(W^{2,1}(\mathbb R^2)\). The regular-value
approximation of \cite{BKK}, in the form of
\cite[Theorem~2.1(iii)]{KPR2015}, supplies \(g_m\in C^1\) and open
\(O_m\subset\mathbb R\) with
\[
 |O_m|<2^{-m},\qquad
 \psi(x)\notin O_m
 \ \Longrightarrow\
 \psi(x)=g_m(x),\quad \nabla\psi(x)=\nabla g_m(x)\ne0
\]
on the patch. Cover \(\{\nabla g_m\ne0\}\) by countably many relatively
compact implicit-function charts. Apply \eqref{eq:5.15} in each such chart.
For \(t\notin O_m\), the portion of \(\{\psi=t\}\) in the chart lies on the
corresponding \(g_m\)-leaf. The Sobolev representatives agree with the
absolutely continuous restrictions on almost every leaf, by slicing.

Enumerate all charts, patches, and integers \(m\), and take a diagonal
subsequence for which the integral in \eqref{eq:5.15} is at most
\(2^{-\ell}\) on each of the first \(\ell\) charts. The exceptional leaf
sets have null union. In addition,
\[
 \left|\limsup_{m\to\infty}O_m\right|
 \le\inf_M\sum_{m\ge M}|O_m|=0.
\]
Outside the resulting null set, a compact regular cycle is covered by
finitely many good charts, and the desired uniform convergence follows.
In space the comparison
\[
 2\pi r_*\int_K|Z|^q\,dr\,dz
 \le\int_{\widetilde K}|Z|^q\,dx
 \le2\pi r^*\int_K|Z|^q\,dr\,dz,
 \qquad 0<r_*\le r\le r^*,
\]
and the analogous derivative bounds reduce the assertion to the same
two-dimensional calculation.
\end{proof}

\subsection{Head values on arcs of the Kronrod tree}

In the plane take the level-component tree \(T_\psi\) after filling the
holes by their constant stream values; see
\cite[Remark~2.4]{KPR2015}. For a meridian separated from the axis, enclose
its closure in a smooth disk in \(r>0\) and use the extended fields. For a
meridian meeting the axis use the compact truncated meridian constructed
below, filling its inner holes by constants. Write \(B_j\) for the element
containing the physical wall, its meridional generator, or its truncated
generator as appropriate. Bernoulli's law assigns a single head value
\(\Phi(C)\) to a positive-diameter level component, through its precise
representative outside an exceptional set of \(\mathcal H^1\)-measure zero.
On the separating arcs used below the following properties hold:
\begin{equation}\label{eq:5.16}
\begin{gathered}
 \Phi(B_j)=p_j,\qquad
 \Phi:[A,B]\longrightarrow\mathbb R
 \ \hbox{is continuous if }\operatorname{diam}A,\operatorname{diam}B>0,\\
 \mathcal H^1(\{\psi(C):C\in Z\})=0
 \ \Longrightarrow\
 \mathcal H^1(\{\Phi(C):C\in Z\})=0,\qquad Z\subset[A,B],\\
 \mathcal H^1(\{\Phi(C):C\in[A,B],\ C\text{ is not regular}\})=0.
\end{gathered}
\end{equation}
Here \emph{regular} includes the convergence in Lemma~\ref{lemma:5.2}.
For the plane, the Euler tree assertions are
\cite[Lemmas~3.5--3.6 and Corollary~3.2]{KPR2015}; for space see
\cite[Lemmas~4.5--4.7 and Corollary~4.1]{KPR2015}.

We spell out the null-set step on interior separating arcs in the presence
of swirl. Let \(G\subset\mathbb R\) be a Borel null set of stream values,
enlarged to include all physical wall values. For \(E=\psi^{-1}(G)\), the
coarea formula gives \(\nabla\psi=0\) a.e.\ on \(E\). On each off-axis patch
of the continued Euler field,
\[
 v_r=v_z=0,\qquad \partial_z\Phi=0
 \quad\hbox{a.e. on }E.
\]
The second assertion is \cite[(3.19) and Remark~3.6]{KPR2015axis}; it does
not discard the azimuthal velocity. Let \(A_v\) be the \(\mathcal H^1\)-null
exceptional set in Bernoulli's law. On almost every line \(r=\rho\), the
precise head agrees outside \(A_v\) with the absolutely continuous
representative \cite[Lemma~3.8]{KPR2015axis}. The one-dimensional area
formula therefore implies
\[
 \mathcal H^1\bigl(\{\Phi(\rho,z):(\rho,z)\in E\setminus A_v\}\bigr)=0
 \quad\hbox{for a.e. }\rho.
\]
Choose a countable dense set of such radii, avoiding also the projection
of \(A_v\), using a countable exhaustion by off-axis rectangles. The union
of the corresponding head images is null.

Every interior element separating the endpoints of a tree arc, and
disjoint from the boundary of the filled planar disk, separates that disk.
Its radial projection is a nondegenerate interval: a compact connected
subset of a vertical line is a point or a segment and cannot separate the
disk. Thus it meets a selected good line. If its stream value is in \(G\),
its head belongs to the preceding null union. An element meeting a
constant filled hole or wall has the corresponding wall head and
contributes at most one of finitely many values. The artificial axial cut
is removed from the selected arcs below. This proves the null-set property
needed in \eqref{eq:5.16}.

Include in \(G\) the exceptional stream values of Lemma~\ref{lemma:5.2}
and of the Sobolev Morse--Sard theorem. The selected regular cycles then
have stream values different from all wall values and lie compactly inside
the physical fluid. The continued Euler geometry thus requires viscous
convergence only where the viscous fields are defined.

\subsection{Geometry at the symmetry axis}

\begin{lemma}[Truncation of an axisymmetric Lipschitz domain]\label{lemma:5.axis}
Suppose that \(\Omega\subset\mathbb R^3\) is bounded, connected,
axisymmetric, and Lipschitz, and meets the axis. There is \(r_0>0\) such that
for every \(0<\varepsilon\le r_0\),
\begin{equation}\label{eq:5.17}
 D_\varepsilon=D\cap\{r>\varepsilon\}
 \quad\hbox{is connected and Lipschitz}.
\end{equation}
Each physical boundary generator meeting the axis remains connected after
truncation and has maximal radius greater than \(2r_0\). Every off-axis
boundary component is contained in \(\{r>2r_0\}\).
\end{lemma}

\begin{proof}
At an axis boundary point choose a local Lipschitz graph direction \(e\)
and a coherent orientation of its normal such that
\(n(y)\cdot e\ge c>0\) a.e.\ in a smaller boundary neighborhood. Use a
ball centered on the axis. Rotational invariance and Fubini give
\begin{equation}\label{eq:5.18}
 n(y)\cdot R_\theta e\ge c\quad\hbox{for a.e. }(y,\theta),
 \qquad
 n(y)\cdot\left(\frac1{2\pi}\int_0^{2\pi}R_\theta e\,d\theta\right)
 =n_z(y)e_z\ge c.
\end{equation}
Hence \(e_z\ne0\), and the axial direction is uniformly transverse.
For completeness, in the original graph coordinates
\(y_n=g(y')\), write the appropriately oriented axial vector as
\((b',b_n)\). The normal estimate yields
\(b_n-b'\cdot\nabla g\ge c_0>0\) a.e.
Convolution of \(g\) on a smaller graph cylinder preserves this inequality.
Integrating on parallel line segments and passing to the uniform limit
shows that \(y_n-g(y')\) increases by at least \(c_0t\) after displacement
by \(t(b',b_n)\). Thus every nearby axial line meets the graph exactly
once. The Lipschitz bound and the lower slope bound give a Lipschitz
intersection height in the horizontal variables. Rotation invariance
reduces it to
\begin{equation}\label{eq:5.19}
 z=\varphi(r),\qquad
 |\varphi(r)-\varphi(s)|\le L|r-s|.
\end{equation}
There is a single axis contact in each such neighborhood. Compactness
therefore gives finitely many contact heights
\(\zeta_1<\cdots<\zeta_m\).

Choose disjoint graph neighborhoods and shrink \(r_0\) so that every
boundary point with \(r\le2r_0\) belongs to one of them. Their graph
functions satisfy
\begin{equation}\label{eq:5.20}
 \varphi_\ell(0)=\zeta_\ell,\qquad
 \varphi_{\ell+1}(r)-\varphi_\ell(r)
 \ge\tfrac12(\zeta_{\ell+1}-\zeta_\ell)>0
 \quad(0\le r\le2r_0).
\end{equation}
There are increasing piecewise affine maps \(H_r\) in the \(z\)-variable
with \(H_r(\zeta_\ell)=\varphi_\ell(r)\) and
\(c\le\partial_zH_r\le C\), such that
\((r,z)\mapsto(r,H_r(z))\) is bi-Lipschitz on the strip. The fluid intervals
correspond to a fixed finite union \(I\) of intervals between consecutive
contact heights:
\[
 D\cap\{0<r<2r_0\}
 =\{(r,H_r(s)):0<r<2r_0,\ s\in I\}.
\]
The meridian \(D\) is connected: paths in the open spatial domain can be
perturbed off the axis and then projected to \((r,z)\).
Given two points of \(D_\varepsilon\), choose a projected path joining
them. Choose \(\eta>0\) with \(\varepsilon+\eta<2r_0\) and with
\(\varepsilon+\eta\) less than both endpoint radii. Inside the strip,
replace its points by
\[
 (r,H_r(s))\longmapsto
 \bigl(\max\{r,\varepsilon+\eta\},
 H_{\max\{r,\varepsilon+\eta\}}(s)\bigr).
\]
The resulting path stays in \(D_\varepsilon\), agrees with the original
path outside the modified strip, and preserves its endpoints.

At an intersection of the cut with a physical graph the local conditions
are \(r>\varepsilon\) and \(\pm(z-\varphi(r))>0\). Along
\(b=(1,\pm(L+1))\), both defining functions increase at least linearly.
Each is a Lipschitz epigraph in this direction, and their intersection is
the epigraph of the maximum of two Lipschitz functions. Elsewhere one has
a single Lipschitz graph. Only the initial segments at the axis are removed
from each axis-meeting generator. Shrinking \(r_0\) below half its positive
maximal radius, and below half the distance of each off-axis component
from the axis, proves the last assertions.
\end{proof}

The local graph representation also proves
\begin{equation}\label{eq:5.21}
 \overline\Omega\cap\{r=0\}\ne\varnothing
 \quad\Longrightarrow\quad
 \Omega\cap\{r=0\}\ne\varnothing.
\end{equation}
Indeed, at an axis boundary point the domain lies on one side of
\eqref{eq:5.19} and contains an adjacent open axial interval.

We use the near-axis Euler assertion
\cite[Lemma~4.4]{KPR2015}: for the compact connected streamfunction
continua in question,
\begin{equation}\label{eq:5.22}
 \min_{C_\ell}r\longrightarrow0,\qquad
 \inf_\ell\max_{C_\ell}r>0
 \quad\Longrightarrow\quad
 \Phi(C_\ell)\longrightarrow p_{\rm ax}.
\end{equation}
For two fixed boundary generators separated by such continua, choose
accessible points on their portions in \(\{r>r_0\}\).
Lemma~\ref{lemma:5.axis}, together with short terminal arcs in boundary
charts, gives a curve \(\gamma\) joining these points, with its interior
in \(D_{r_0}\). Every separating continuum intersects \(\gamma\); hence
\(\sup_Cr\ge\min_\gamma r>0\). For an interior terminal continuum use
one endpoint on that continuum and first decrease \(r_0\).
Therefore, a fixed range of head values strictly above \(p_{\rm ax}\)
cannot contain separating continua approaching the axis.

This is the axis-removal input in
\cite[Lemmas~4.6 and 4.10]{KPR2015}. The selected tree arcs, cycles, and
zero-head wall components may consequently be restricted to
\begin{equation}\label{eq:5.23}
 r\ge r_*>0.
\end{equation}
The cut is fixed before the indices \(i,k\) below are chosen.
The continued Euler fields satisfy the cited near-axis assertion on the
smooth auxiliary domain; see also \cite[Lemma~3.12]{KPR2015axis}.
On the physical truncated meridian, Lemma~\ref{lemma:5.axis} and the
curve \(\gamma\) verify the geometric conditions (4.24)--(4.25) in
\cite{KPR2015}. If the claimed fixed cut failed, the argument of its
Lemmas~4.6 and 4.10 would produce separating elements meeting cuts
\(r\downarrow0\), with uniformly positive radial extent and head tending
to the terminal value, which is strictly above \(p_{\rm ax}\). This
contradicts \eqref{eq:5.22}.

The continuity on tree arcs used here is obtained as in
\cite[Lemmas~3.5 and 4.5]{KPR2015}: intermediate elements separate the
endpoints; Hausdorff limits contain positive-diameter continua in the
limiting element; and an absolutely continuous good line through a common
projection gives convergence of the constant head values. The precise
representatives are those of the continued Euler fields, including across
the original boundary. No boundary regularity of the viscous sequence is
used in this step.

\subsection{The separating-level construction}\label{subsec:separating}

Assume the boundary-maximum alternative, and shift the pressure so that
\begin{equation}\label{eq:5.24}
 \operatorname*{ess\,sup}_\Omega\Phi
 =\max_jp_j=0,\qquad \min_jp_j<0.
\end{equation}
In the spatial boundary-maximum argument assume that the zero-pressure
components are off-axis. If the domain meets the axis, use
\eqref{eq:5.23}.

Fix a terminal tree element \(B_+\) containing a zero-pressure boundary
component, so that \(\Phi(B_+)=0\), and connect each negative-pressure
component to \(B_+\). After axis removal, the initial elements can instead
be the finitely many negative-head elements supplied by
\cite[Lemma~4.6]{KPR2015}. Write the resulting arcs as
\(\mathcal A_j=[C_j,B_+]\), with \(\Phi(C_j)<0\). Let
\(\mathcal N\subset\mathbb R\) be the null set of exceptional head values
on their union. Choose
\begin{equation}\label{eq:5.25}
 0<t_0<\min_j(-\Phi(C_j)),\qquad
 t_0\notin\bigcup_{i=0}^{\infty}2^i(-\mathcal N),\qquad
 t_i=2^{-i}t_0.
\end{equation}
This choice is possible because the forbidden union is null. In the order
toward \(B_+\), set
\begin{equation}\label{eq:5.26}
 A_i^j=\max\{C\in\mathcal A_j:\Phi(C)=-t_i\}.
\end{equation}
Continuity on the arc gives
\[
 \Phi(A_i^j)=-t_i,\qquad
 \Phi(C)>-t_i\quad(C\in(A_i^j,B_+]),\qquad
 A_i^j<A_{i+1}^j<B_+.
\]
All these elements are regular interior cycles of the physical
streamfunction: the excluded head values include the images of the
finitely many stream values of the added constant regions. Repeated cycles
are counted only once.

Let \(V_i\) be the connected region on the \(B_+\)-side of the selected
cycles. Jordan separation and the tree order imply
\begin{equation}\label{eq:5.27}
 V_{i+1}\subset V_i,\qquad
 \partial V_i=\Gamma(i)\cup\bigcup_jA_i^j,
\end{equation}
where \(\Gamma(i)\) consists only of zero-pressure physical components.
Indeed, a negative wall is separated from \(B_+\) by one selected cycle;
a sufficiently small fluid-side collar of the terminal wall lies on the
opposite side. Every other wall lies wholly on one side of each cycle,
by connectedness. Filling a hole for the tree construction does not join
different fluid components on a selected side: each path segment through
the hole may be replaced by a path in a connected fluid-side collar of its
boundary, disjoint from the finitely many cycles. Equivalently, the regions
are obtained by successive Jordan cuts in the physical domain.

In a truncated meridian the discarded axial strip lies on the outer side
of the cut separating the terminal component from the axial boundary;
compare \cite[(4.26)--(4.28)]{KPR2015}. If the terminal set is an interior
continuum, take cuts on the arcs from every physical wall. Then no physical
wall remains on the terminal side, giving the compact interior regions
used later.

The decreasing sequence of subsets of the finite set of wall components
stabilizes. After discarding finitely many indices, write
\begin{equation}\label{eq:5.28}
 \Gamma(i)=\Gamma=\bigcup_{j\in J_0}\Gamma_j,\qquad
 J_0\ne\varnothing,\qquad p_j=0\quad(j\in J_0).
\end{equation}
These are the constructions of
\cite[(3.29)--(3.31) and (4.28)]{KPR2015}. In spatial formulas the cycles
and planar regions are identified with their rotations, unless meridional
measure is explicitly indicated.

Use the open separating layer
\begin{equation}\label{eq:5.29}
 L_i=V_i\setminus\overline{V_{i+1}},\qquad
 I_i=[5t_i/8,\,7t_i/8],\qquad |I_i|=t_i/4.
\end{equation}
For a fixed \(i\), both regions contain a sufficiently small collar of their
common physical boundary \(\Gamma\). Every other boundary is a finite union
of interior cycles. Hence
\(\overline{L_i}\Subset\Omega\), and in space
\(L_i\subset\{r\ge r_*\}\). The cycle boundaries have zero volume, so
\begin{equation}\label{eq:5.30}
 |L_i|=|V_i|-|V_{i+1}|\longrightarrow0.
\end{equation}
No assertion that \(|\bigcap_iV_i|=0\) is needed.

By Lemma~\ref{lemma:5.2}, for \(k\ge k_i\),
\[
 \max_j\|\Phi_k+t_i\|_{C(A_i^j)}<t_i/16,\qquad
 \max_j\|\Phi_k+t_{i+1}\|_{C(A_{i+1}^j)}<t_i/16.
\]
In particular,
\begin{equation}\label{eq:5.31}
 \Phi_k<-7t_i/8\quad\hbox{on }\bigcup_jA_i^j,\qquad
 \Phi_k>-5t_i/8\quad\hbox{on }\bigcup_jA_{i+1}^j.
\end{equation}
Choose a fixed open anchor in \(L_i\), adjacent to all inner cycles and
contained in \(\{\Phi_k>-5t_i/8\}\). For \(t\in I_i\), let \(W_i^k(t)\)
be the union of the components of \(\{\Phi_k>-t\}\cap L_i\) meeting
this anchor, and set \(S_i^k(t)=\partial_{L_i}W_i^k(t)\).

For almost every \(t\), the Sobolev Morse--Sard theorem applied to the
individual interior \(H^2\) head gives finitely many regular closed curves,
or their rotated tori. The outer normal of the chosen superlevel region
points toward decreasing head. The gradient in the flux integral is its
Sobolev trace from the interior. This agrees surface-a.e.\ with the
precise representative: locally \(\Phi_k\in H^2\), so \(\nabla\Phi_k\in
H^1\), and its trace on a \(C^1\) curve or an off-axis rotated surface
agrees with those precise values. Thus
\begin{equation}\label{eq:5.32}
 \Phi_k=-t,\qquad
 n=-\frac{\nabla\Phi_k}{|\nabla\Phi_k|},\qquad
 \nabla\Phi_k\cdot n=-|\nabla\Phi_k|
 \quad\hbox{on }S_i^k(t).
\end{equation}

To include the common cycle interfaces explicitly, define the glued open
region by
\begin{equation}\label{eq:5.33}
 D_i^k(t,h)=
 \operatorname{int}_{\Omega}
 \bigl(\overline{W_i^k(t)}\cup\overline{V_{i+1}}\bigr)
 \setminus\overline{C_h}.
\end{equation}
Here \(C_h\) is the collar of \(\Gamma\), and the interior is relative to
\(\Omega\). The strict inner-cycle inequality in \eqref{eq:5.31} supplies
full neighborhoods on both sides of every common cycle; hence no inner
cycle remains as a crack. For sufficiently small \(h>0\),
\begin{equation}\label{eq:5.34}
 D_i^k(t,h)\Subset\Omega,\qquad
 \partial D_i^k(t,h)=S_i^k(t)\,\dot\cup\,\Gamma_h.
\end{equation}
Each component has finitely many regular interior surfaces and a Lipschitz
collar surface as its boundary. Gauss--Green therefore applies componentwise.

The set \(J_0\) is nonempty and proper, since at least one negative-pressure
component has been separated from the terminal wall. Choose
\(d_0,\ldots,d_N\in\mathbb R\) with
\(\sum_{j=0}^Nd_j=0\) and \(\sum_{j\in J_0}d_j=1\).
For \(j\ge1\), choose \(y_j\in\Omega_j\), set
\(s_{n-1}=\mathcal H^{n-1}(\mathbb S^{n-1})\), and define
\[
 Z(x)=-\frac1{s_{n-1}}\sum_{j=1}^N
          d_j\,\frac{x-y_j}{|x-y_j|^n}.
\]
This field is smooth, bounded, and solenoidal near \(\overline\Omega\).
Its fluid-oriented flux through \(\Gamma_j\) is \(d_j\): integrate in
the \(j\)-th hole after removing a small ball about \(y_j\); the outer
flux follows by incompressibility. Applying Gauss--Green first in the
removed collars and then in \(D_i^k(t,h)\) gives
\[
 \int_{\Gamma_h}Z\cdot n_D\,dS=1,\qquad
 \int_{S_i^k(t)}Z\cdot n_D\,dS=-1.
\]
Consequently, in both dimensions,
\begin{equation}\label{eq:5.35}
 \mathcal H^{n-1}(S_i^k(t))
 \ge c_*:=\|Z\|_{L^\infty(\Omega)}^{-1}>0.
\end{equation}

\subsection{An interior lower bound for the enstrophy}

For each fixed \(i\), the Euler head is nonconstant in \(V_{i+1}\).
Indeed, the regular cycles at levels \(-t_{i+2}\) and \(-t_{i+3}\)
lie inside this region and have different head values. If \(\omega=0\)
a.e.\ in the connected region \(V_{i+1}\), then \eqref{eq:5.4} would imply
\(\nabla\Phi=0\) a.e.\ and \(\Phi=c\) a.e.\ there. The precise representative
of an a.e.\ constant on an open set is that same constant; choosing on each
cycle a point outside the \(\mathcal H^1\)-null Bernoulli exceptional set
would contradict the two distinct values. Therefore
\(\int_{V_{i+1}}|\omega|^2\,dx>0\).

Choose a compact \(K_i\subset V_{i+1}\) and \(c_i>0\) with
\begin{equation}\label{eq:5.36}
 \int_{K_i}|\omega|^2\,dx=4c_i>0.
\end{equation}
Weak convergence of the curls gives
\(\liminf_k\int_{K_i}|\omega_k|^2\,dx\ge4c_i\).
Choose \(k_i'\) such that this integral is at least \(2c_i\) for
\(k\ge k_i'\), and choose \(\delta_i>0\) with
\(K_i\cap\overline{C_{\delta_i}}=\varnothing\).
Uniformly for regular \(t\in I_i\) and \(0<h\le\delta_i\),
\begin{equation}\label{eq:5.37}
 K_i\subset D_i^k(t,h),\qquad
 \int_{D_i^k(t,h)}|\omega_k|^2\,dx\ge2c_i\ge c_i
 \quad(k\ge k_i').
\end{equation}

\subsection{Flux identities and the parameter order}

All normals in this subsection point out of \(D=D_i^k(t,h)\).
On \(\Gamma_h\), they point toward the removed physical collar.
Gauss--Green in that collar gives
\begin{equation}\label{eq:5.38}
 \int_{\Gamma_h}u_k\cdot n_D\,dS
 =\int_\Gamma\operatorname{Tr}u_k\cdot n_\Omega\,dS
 =\frac{\lambda_k\mu_k}{\nu}\sum_{j\in J_0}F_j.
\end{equation}
Since \(\operatorname{div}u_k=0\) in \(D\),
\begin{equation}\label{eq:5.39}
\begin{gathered}
 \int_{S_i^k(t)}u_k\cdot n_D\,dS
 =-\frac{\lambda_k\mu_k}{\nu}\sum_{j\in J_0}F_j,\\
 \mu_k^{-1}
 \left|\int_{S_i^k(t)}u_k\cdot n_D\,dS\right|
 \le C_F,\qquad
 C_F=\nu^{-1}\sum_{j=0}^N|F_j|.
\end{gathered}
\end{equation}
The constant is independent of every selection of components and indices.

Fix \(i\) first. Let \(C_N\) be the coefficient of the velocity product
in \eqref{eq:4.27}. Select, in this order,
\begin{equation}\label{eq:5.40}
\begin{gathered}
 0<\varepsilon_i<\frac{c_i}{4(C_N+1)},\qquad
 C_{\varepsilon_i}\ \hbox{from Lemma~\ref{lemma:4.3}},\\
 0<\sigma_i<\frac{c_i}{4\sqrt{C_{\varepsilon_i}}+1},\qquad
 0<\delta\le\min\{\delta_i,\delta_{\sigma_i},h_0/2\}.
\end{gathered}
\end{equation}
Further decrease \(\delta\) to avoid the finitely many separating cycles.
Then take \(k\) large enough for all earlier requirements and so that the
force errors below total at most \(c_i/4\).
On the common selected surface \(h_k\in[\mu_k^2,\delta]\),
\begin{equation}\label{eq:5.41}
\begin{gathered}
 \|\Phi_k\|_{L^2(\Gamma_{h_k})}\le\sigma_i,\qquad
 \|u_k\|_{L^2(\Gamma_{h_k})}
 \le\sqrt{C_{\varepsilon_i}}\,\mu_k,\\
 \int_{\Gamma_{h_k}}|u_k||\nabla u_k|\,dS\le\varepsilon_i.
\end{gathered}
\end{equation}
Thus
\begin{equation}\label{eq:5.42}
\begin{aligned}
 \mu_k^{-1}\left|\int_{\Gamma_{h_k}}\Phi_ku_k\cdot n_D\,dS\right|
 &\le\sigma_i\sqrt{C_{\varepsilon_i}}\le c_i/4,\\
 \int_{\Gamma_{h_k}}\nabla\Phi_k\cdot n_D\,dS
 &\le C_N\varepsilon_i+C\mu_k^2\le c_i/4+C\mu_k^2.
\end{aligned}
\end{equation}
For every measurable \(D\subset\Omega\),
\begin{equation}\label{eq:5.43}
\begin{aligned}
 \mu_k^{-1}\left|\int_Dg_k\cdot u_k\,dx\right|
 &\le\mu_k^{-1}\|g_k\|_2\|u_k\|_2\le C\mu_k,\\
 \left|\int_D\operatorname{div}g_k\,dx\right|
 &\le|\Omega|^{1/2}\|\operatorname{div}g_k\|_2\le C\mu_k^2.
\end{aligned}
\end{equation}
The complete order is
\[
 i\ \longrightarrow\ c_i\ \longrightarrow\ \varepsilon_i
 \ \longrightarrow\ C_{\varepsilon_i}\ \longrightarrow\ \sigma_i
 \ \longrightarrow\ \delta\ \longrightarrow\ k\ \longrightarrow\ h_k.
\]
In particular, \(\delta\) is fixed as \(k\to\infty\) in the
surface-selection lemma, and \(h_k\) is independent of \(t\in I_i\).

\subsection{The head integral and the coarea contradiction}

Write \(S=S_i^k(t)\) and \(D=D_i^k(t,h_k)\).
Integration of \eqref{eq:4.26} gives the signed identity
\begin{equation}\label{eq:5.44}
\begin{aligned}
 -\int_S|\nabla\Phi_k|\,dS
 +\int_{\Gamma_{h_k}}\nabla\Phi_k\cdot n_D\,dS
 &=\int_D|\omega_k|^2\,dx
   -\frac{t}{\mu_k}\int_Su_k\cdot n_D\,dS\\
 &\quad+\frac1{\mu_k}
       \int_{\Gamma_{h_k}}\Phi_ku_k\cdot n_D\,dS\\
 &\quad-\frac1{\mu_k}\int_Dg_k\cdot u_k\,dx
       +\int_D\operatorname{div}g_k\,dx.
\end{aligned}
\end{equation}
Rearranging and using \eqref{eq:5.37}--\eqref{eq:5.43}, we obtain
\begin{equation}\label{eq:5.45}
\begin{aligned}
 \int_{S_i^k(t)}|\nabla\Phi_k|\,dS
 &\le C_Ft-c_i+C_N\varepsilon_i
       +\sigma_i\sqrt{C_{\varepsilon_i}}+C\mu_k+C\mu_k^2\\
 &\le C_Ft.
\end{aligned}
\end{equation}
This holds for a.e.\ \(t\in I_i\), with a constant independent of \(i\).
The fixed positive enstrophy absorbs the errors before integration in \(t\).

We record measurability of the swept set. Fix \(i,k\), and let
\(C\subset L_i\) be the open anchor chosen above; it satisfies
\[
 C\subset\{\Phi_k>-5t_i/8\}\subset\{\Phi_k>-t\}
 \quad(t\in I_i).
\]
Define
\(\mathcal A=\{(x,t)\in L_i\times I_i:x\in W_i^k(t)\}\).
If \((x,t)\in\mathcal A\), a compact path \(\gamma\) in the open superlevel
set joins \(x\) to \(C\), and
\(\min_{y\in\gamma}(\Phi_k(y)+t)>0\).
The strict margin persists for nearby \(x,t\), so \(\mathcal A\) is
relatively open. Let
\[
 d(x,t)=\operatorname{dist}(x,W_i^k(t)),\qquad
 \mathcal B=\{(x,t):d(x,t)=0\}\setminus\mathcal A.
\]
The distance is upper semicontinuous: when \(d(x,t)<a\), a point of
\(W_i^k(t)\) within distance \(a\) remains in nearby sections by the
openness of \(\mathcal A\). Thus \(\mathcal B\) is Borel, with sections
\(\mathcal B_t=S_i^k(t)\). Choose a Borel full-measure set
\(R_{ik}\subset I_i\) on which regularity and all level estimates hold.
Then
\begin{equation}\label{eq:5.46}
 E_{ik}=\{x\in L_i:-\Phi_k(x)\in R_{ik},\
                  (x,-\Phi_k(x))\in\mathcal B\}
\end{equation}
is Borel, and
\(E_{ik}\cap\{\Phi_k=-t\}=S_i^k(t)\) for \(t\in R_{ik}\).
In space perform the construction meridionally and pull it back under the
continuous map \(x\mapsto(r,z)\).

Coarea and \eqref{eq:5.45} give
\begin{equation}\label{eq:5.47}
\begin{aligned}
 \int_{E_{ik}}|\nabla\Phi_k|^2\,dx
 &=\int_{R_{ik}}\int_{S_i^k(t)}|\nabla\Phi_k|\,dS\,dt\\
 &\le C_F\int_{5t_i/8}^{7t_i/8}t\,dt
 =\frac{3C_F}{16}t_i^2.
\end{aligned}
\end{equation}
The lower surface bound \eqref{eq:5.35} also gives
\begin{equation}\label{eq:5.48}
 \int_{E_{ik}}|\nabla\Phi_k|\,dx
 =\int_{R_{ik}}\mathcal H^{n-1}(S_i^k(t))\,dt
 \ge\frac{c_*t_i}{4}.
\end{equation}
With \(C_A=1+3C_F/16>0\), Cauchy--Schwarz and \(E_{ik}\subset L_i\)
imply
\begin{equation}\label{eq:5.49}
 \frac{c_*t_i}{4}
 \le\left(\int_{E_{ik}}|\nabla\Phi_k|^2\,dx\right)^{1/2}
          |E_{ik}|^{1/2}
 \le\sqrt{C_A}\,t_i|L_i|^{1/2},
 \qquad
 |L_i|\ge\frac{c_*^2}{16C_A}>0.
\end{equation}
This contradicts \eqref{eq:5.30}. For each \(i\), one first chooses
\(k\ge k(i)\); the lower bound is independent of this choice.

\subsection{An interior head value above all boundary pressures}

Suppose that
\begin{equation}\label{eq:5.50}
 \operatorname*{ess\,sup}_\Omega\Phi>\max_jp_j,
\end{equation}
allowing an infinite essential supremum. The interior-continuum
construction of \cite[Lemma~3.10 and Section~4.2.3]{KPR2015} supplies a
compact connected interior set \(F\) of positive diameter on a
streamfunction level, with \(\Phi(F)>\max_jp_j\).

To explain localization, put \(b=\max_jp_j\). A positive-measure set of
physical points satisfies \(\Phi>b+\eta\) for some \(\eta>0\).
Choose a good vertical line meeting this set, with an absolutely continuous
head on compact segments and avoiding the Bernoulli exceptional set; in
space take \(r=\rho>0\). Outside the enclosing domain the continued head
is \(p_0\le b\). Continuity gives a compact subsegment on which
\(\Phi>b\) and on which \(\Phi\) takes at least two values. This segment
lies in the fluid, because the good-line values on walls and complementary
regions are their wall pressures. Bernoulli's law forces \(\psi\) to be
nonconstant there. Its stream values contain an interval; choose an
interior value outside the Morse--Sard exception and the finite set of
wall values. The regular level through a corresponding point contains a
small compact arc with constant head above \(b\). This supplies \(F\) in
the meridian without requiring an entire level component to have already
been separated from the axis. In the plane the corresponding full
component is compactly inside the fluid. For a spatial domain meeting the
axis, apply \cite[Lemma~4.10]{KPR2015} using
Lemma~\ref{lemma:5.axis}. Shift by \(\Phi(F)\), so that
\(\Phi(F)=0\) and \(\max_jp_j<0\).

In the plane \(F\) can be taken to be a tree element. On the truncated
meridian use the element containing \(F\) as the common endpoint of the
arcs. The last-crossing construction gives
\begin{equation}\label{eq:5.51}
 V_i\Subset\Omega,\qquad
 \partial V_i=\bigcup_jA_i^j,\qquad
 \Phi|_{A_i^j}=-t_i.
\end{equation}
There are no physical boundary components. For fixed \(i\), sufficiently
large \(k\), and regular \(t\in I_i\), define
\[
 D_i^k(t)=\operatorname{int}_{\Omega}
       \bigl(\overline{W_i^k(t)}\cup\overline{V_{i+1}}\bigr).
\]
Then \(D_i^k(t)\Subset\Omega\),
\(\partial D_i^k(t)=S_i^k(t)\), and \(\Phi_k=-t\) on the entire boundary.
Two later regular cycles show that \(\Phi\) is nonconstant in \(V_{i+1}\),
as before. A compact \(K_i\subset V_{i+1}\) therefore gives the lower
enstrophy bound \(c_i>0\) for all large \(k\).
Incompressibility cancels the complete convective boundary term:
\begin{equation}\label{eq:5.52}
 \int_{\partial D_i^k(t)}\Phi_ku_k\cdot n\,dS
 =-t\int_{\partial D_i^k(t)}u_k\cdot n\,dS
 =-t\int_{D_i^k(t)}\operatorname{div}u_k\,dx=0.
\end{equation}
Integration of the head equation now yields
\begin{equation}\label{eq:5.53}
\begin{aligned}
 0&>-\int_{S_i^k(t)}|\nabla\Phi_k|\,dS\\
 &=\int_{D_i^k(t)}|\omega_k|^2\,dx
   -\mu_k^{-1}\int_{D_i^k(t)}g_k\cdot u_k\,dx
   +\int_{D_i^k(t)}\operatorname{div}g_k\,dx\\
 &\ge c_i-C\mu_k-C\mu_k^2>0
 \quad\hbox{for sufficiently large }k.
\end{aligned}
\end{equation}
No physical-boundary estimate is needed in this alternative.

\subsection{The axisymmetric boundary maximum at the axis}

It remains to consider a spatial domain meeting the axis with
\(\operatorname*{ess\,sup}_\Omega\Phi=\max_jp_j=p_{\rm ax}\).
The exterior component meets the axis, so \(p_0=p_{\rm ax}\).
After a pressure shift,
\begin{equation}\label{eq:5.54}
 \Phi\le0\quad\hbox{a.e.},\qquad p_0=0,\qquad
 p_j\le0,\qquad p_{j_*}<0\quad\hbox{for some }j_*\ge1.
\end{equation}
Let \(X=xp+(v\cdot x)v\). The product bounds in \eqref{eq:5.8} give
\(xp\in W^{1,3/2}\), while \(v\in L^6\) and \(\nabla v\in L^2\)
give \((v\cdot x)v\in W^{1,3/2}\). Compute
\begin{equation}\label{eq:5.55}
\begin{aligned}
 \operatorname{div}(xp)&=3p+x\cdot\nabla p,\\
 \operatorname{div}((v\cdot x)v)
 &=\nabla(v\cdot x)\cdot v+(v\cdot x)\operatorname{div}v\\
 &=x\cdot(v\cdot\nabla)v+|v|^2,\\
 \operatorname{div}X&=3p+|v|^2
 =3\Phi-\tfrac12|v|^2.
\end{aligned}
\end{equation}
Since the velocity trace vanishes, Gauss--Green gives
\begin{equation}\label{eq:5.56}
 \int_\Omega(3\Phi-\tfrac12|v|^2)\,dx
 =\int_{\partial\Omega}\operatorname{Tr}X\cdot n_\Omega\,dS
 =\sum_{j=0}^Np_j\int_{\Gamma_j}x\cdot n_\Omega\,dS.
\end{equation}
On an inner boundary, \(n_\Omega=-n_{\Omega_j}\), whence
\(\int_{\Gamma_j}x\cdot n_\Omega\,dS=-3|\Omega_j|\).
Because \(p_0=0\), \eqref{eq:5.54} implies
\begin{equation}\label{eq:5.57}
 0\ge\int_\Omega(3\Phi-\tfrac12|v|^2)\,dx
 =-3\sum_{j=1}^Np_j|\Omega_j|>0.
\end{equation}
This excludes the remaining alternative. If \(\max_jp_j>p_{\rm ax}\),
all maximal-pressure components avoid the axis and the preceding
boundary-maximum proof applies. If \(\Omega\) does not meet the axis,
\eqref{eq:5.21} implies \(\min_{\overline\Omega}r>0\), so all meridional
weights are globally comparable and no axial truncation is needed.

\begin{proposition}[Exclusion of a normalized blow-up sequence]\label{prop:5.4}
Under the normalization, energy identity, and pressure and collar estimates
established in the preceding sections, no blow-up sequence
\eqref{eq:5.1} exists in the planar or axisymmetric spatial setting.
\end{proposition}

\begin{proof}
Let \(M=\operatorname*{ess\,sup}_\Omega\Phi\) and \(P=\max_jp_j\).
The trace implies \(M\ge P\): if \(M<P\), the a.e.\ bound
\(\Phi\le M\) would give \(\operatorname{Tr}\Phi\le M\) on every
component, contradicting one with pressure \(P\). The alternatives are
\[
\begin{array}{lll}
 M>P
 &\quad\Longrightarrow&\text{contradiction \eqref{eq:5.53}},\\
 M=P,\ n=2
 &\quad\Longrightarrow&\text{contradiction \eqref{eq:5.49}},\\
 M=P,\ n=3,\ \Omega\cap\{r=0\}=\varnothing
 &\quad\Longrightarrow&\text{contradiction \eqref{eq:5.49}},\\
 M=P>p_{\rm ax},\ n=3
 &\quad\Longrightarrow&\text{contradiction \eqref{eq:5.49}},\\
 M=P=p_{\rm ax},\ n=3
 &\quad\Longrightarrow&\text{contradiction \eqref{eq:5.57}}.
\end{array}
\]
The boundary-maximum alternatives have at least one negative boundary
pressure after normalization, by \eqref{eq:5.6}. These cases exhaust all
possibilities.
\end{proof}

\section{Existence for regular forces}\label{sec:regular-existence}

\begin{proof}[Proof of Theorem~\ref{thm:regular}]
Fix $(\Omega,\nu,a,F)$ as in the statement and use the space $X$
of Section~\ref{sec:continuation}. Let
\[
 \mathcal S=\{w\in X:w=\lambda T(w)\text{ for some }\lambda\in[0,1]\}.
\]
If $\sup_{w\in\mathcal S}\|\nabla w\|_2=\infty$, choose
$w_k\in\mathcal S$ with $\|\nabla w_k\|_2\ge k$.
Equations~\eqref{eq:2.8} and \eqref{eq:2.10} give the normalized
sequence of Section~\ref{sec:pressure}.
Lemma~\ref{lemma:3.5} controls the regular pressure on every
admissible component. Equations~\eqref{eq:4.12}, \eqref{eq:4.19}
and Lemma~\ref{lemma:4.3} give the head and velocity estimates on
one common surface. Proposition~\ref{prop:5.4} excludes such a
sequence. Consequently
\[
 M=M(\Omega,\nu,a,F;X)
 :=\sup_{w\in\mathcal S}\|\nabla w\|_2<\infty.
\]
Take $R=M+1$ and $B_R=\{w\in X:\|\nabla w\|_2<R\}$.
Then $0\notin(I-\lambda T)(\partial B_R)$ for $0\le\lambda\le1$.
Compactness of $T$ and homotopy invariance of Leray--Schauder
degree yield
\[
 \deg(I-T,B_R,0)=\deg(I,B_R,0)=1,
 \qquad \text{hence some }w\in B_R\text{ satisfies }w=T(w).
\]
For $u=U+w$, equations~\eqref{eq:2.3} and \eqref{eq:2.6} give
\[
 \nu(\nabla u,\nabla\eta)+b(u,u,\eta)=\langle F,\eta\rangle
 \quad(\eta\in V),\qquad\operatorname{div}u=0,
 \qquad\operatorname{Tr}u=a.
\]
In the symmetric case the equation on all of $V$ follows from
the group averaging argument in Section~\ref{sec:continuation}.
For $\varphi\in H$, define
\[
 \ell(\varphi)=\nu(\nabla u,\nabla\varphi)
 -(u\otimes u,\nabla\varphi)-\langle F,\varphi\rangle.
\]
It vanishes on $V$, and
\[
 |\ell(\varphi)|\le
 (\nu\|\nabla u\|_2+\|u\|_4^2+\|F\|_{H^{-1}})
 \|\nabla\varphi\|_2.
\]
Formula~\eqref{eq:2.4} gives the unique mean-zero pressure
$P_F\in L_0^2$ satisfying \eqref{eq:1.full}. In space, $u$ is
fixed by the prescribed group; since $F$ is also invariant,
the mean-zero pressure has the same invariance by uniqueness.
\end{proof}

The bound $M$ is for fixed data. The spatial argument uses
subcritical trace compactness on the generating curve of an
off-axis boundary component. For a general two-dimensional
boundary surface, by contrast, the indices give
\[
 \frac1{3/2}-\frac{1/3}{2}=\frac12,
 \qquad W^{1/3,3/2}(\partial\Omega)\hookrightarrow L^2(\partial\Omega),
\]
which is critical. Both the compactness argument and the
level-set argument use the meridional structure in the spatial
case. No existence assertion for arbitrary three-dimensional
geometry follows from this proof.

\section{The solenoidal-dual force and the final limit}
\label{sec:energyforce}\label{sec:dual}

Fix the domain, viscosity and trace of Theorem~\ref{thm:main}, and
let $f\in V'$ have the prescribed symmetries. Throughout this
section fix the single representative
\[
 F=\mathcal Ef\in H^{-1},\qquad
 \|F\|_{H^{-1}}=\|f\|_{V'}.
\]
The pressures $P_k$ and $\beta_k$ below belong respectively to
the solution with smooth representative $F_k$ and to a small
zero-boundary correction. Their scaled difference $p_k$ is
the pressure for the corrected velocity $w_k$.
All three are initially normalized to have mean zero.
Later simultaneous shifts of the normalized heads are used
only inside the contradiction argument.

\subsection{Smooth approximation and the normalized energy identity}

Let $z_f\in V$ be given by Lemma~\ref{lem:extension}. Choose
$z_k\in C_c^\infty(\Omega;\mathbb R^n)$ with $z_k\to z_f$ in $H$.
No solenoidal constraint on these smooth approximants is needed.
Since $F=-\Delta z_f$,
\[
 \|{-\Delta z_k}-F\|_{H^{-1}}
 \le\|\nabla(z_k-z_f)\|_2\to0.
\]
In the plane set $F_k=-\Delta z_k$.
In space average over the compact group $G$:
\[
 F_k=\int_G\rho_R(-\Delta z_k)\,dR.
\]
The union of the rotated supports is compactly contained in
$\Omega$, and the action is isometric on $H^{-1}$. Thus in
both dimensions we have symmetric representatives
$F_k\in C_c^\infty(\Omega;\mathbb R^n)$ with
$F_k\to F$ strongly in $H^{-1}$.
Here and below $F_k$ denotes these unscaled representatives,
not the normalized Stokes right-hand sides of \eqref{eq:2.13}.
Theorem~\ref{thm:regular} gives solutions $(U_k,P_k)$ with
force $F_k$, trace $a$, mean-zero pressure and the required
velocity symmetries.

Use the homogeneous unit-viscosity Stokes lift $(A,\pi)$ from
\eqref{eq:2.12}:
\[
\begin{gathered}
 -\Delta A+\nabla\pi=0,\qquad\operatorname{div}A=0,
 \qquad\operatorname{Tr}A=a,\\
 \int_\Omega\pi=0,
 \qquad A\in H^1,\quad\pi\in L^2,\quad\Delta\pi=0.
\end{gathered}
\]
Suppose, after extraction, that
$J_k=\|\nabla(U_k-A)\|_2\to\infty$. Put
$v_k=(U_k-A)/J_k$ and $\mu_k=\nu/J_k$. Along a subsequence,
$v_k\rightharpoonup v$ in $V$ and $v_k\to v$ in $L^4$.
Testing the equation for $U_k$ by $U_k-A$ gives exactly
\begin{equation}\label{eq:7.1}
\begin{split}
 \nu={}&-b(v_k,A,v_k)-J_k^{-1}b(A,A,v_k)\\
 &+J_k^{-1}\bigl(\langle F_k,v_k\rangle
                  -\nu(\nabla A,\nabla v_k)\bigr)
 \longrightarrow b(v,v,A).
\end{split}
\end{equation}
Here $b(v,A,v)=-b(v,v,A)$ and
$|\langle F_k,v_k\rangle|\le\sup_k\|F_k|_V\|_{V'}
\le\sup_k\|F_k\|_{H^{-1}}$. In particular $v\ne0$.

\subsection{The nonlinear correction and its effective force}

Relabel the selected solution--force pairs together, so the
relabelled $F_k$ still converge to $F$. Set
$M_m=\max_{1\le j\le m}\|F_j\|_\infty$ and, for large $k$, let
\begin{equation}\label{eq:7.2}
 m(k)=\max\{m\le k:M_m\le\sqrt{J_k}\},
 \qquad s_k=F_{m(k)},\qquad \rho_k=F_k-s_k.
\end{equation}
The set is eventually nonempty. For each fixed $m$, both
$m\le k$ and $M_m\le\sqrt{J_k}$ eventually hold, hence
$m(k)\to\infty$ and
\[
 J_k^{-1}\|s_k\|_\infty\le J_k^{-1/2}\to0,
 \qquad
 \|\rho_k\|_{H^{-1}}
 \le\|F_k-F\|_{H^{-1}}+\|F-F_{m(k)}\|_{H^{-1}}\to0.
\]
Solve the zero-boundary problem
\begin{equation}\label{eq:7.3}
 -\nu\Delta B_k+(B_k\cdot\nabla)B_k+\nabla\beta_k=\rho_k,
 \qquad B_k\in V.
\end{equation}
Indeed, the variational Stokes inverse applied to
$\rho_k-\operatorname{div}(z\otimes z)$ defines a map $T_k$ on
$V$ satisfying
\[
 \|T_kz\|_V\le\nu^{-1}\|\rho_k|_V\|_{V'}
                 +C\nu^{-1}\|z\|_V^2,
\]
\[
 \|T_kz-T_ky\|_V
 \le C\nu^{-1}(\|z\|_V+\|y\|_V)\|z-y\|_V.
\]
For large $k$ it is a contraction on the closed ball of radius
$2\nu^{-1}\|\rho_k|_V\|_{V'}$ (with the zero-force case immediate).
Its unique small fixed point inherits the data symmetries.
Testing with $B_k$ improves the bound to
\begin{equation}\label{eq:7.4}
 \nu\|\nabla B_k\|_2^2=\langle\rho_k,B_k\rangle,
 \qquad
 \|\nabla B_k\|_2\le\nu^{-1}\|\rho_k|_V\|_{V'}
 \le\nu^{-1}\|\rho_k\|_{H^{-1}}\to0.
\end{equation}
Formula~\eqref{eq:2.4} recovers $\beta_k\in L_0^2$.
The bootstrap in \eqref{eq:4.24}--\eqref{eq:4.25}, followed by
iteration with the smooth forces, gives
\[
 U_k,B_k\in C^\infty_{\rm loc}(\Omega),
 \qquad P_k,\beta_k\in C^\infty_{\rm loc}(\Omega)
 \quad\text{for each fixed }k.
\]
Set $w_k=(U_k-B_k)/J_k$ and $p_k=(P_k-\beta_k)/J_k^2$.
The cancellation of the quadratic correction reads
\[
 (U_k\cdot\nabla)U_k-(B_k\cdot\nabla)B_k
 =J_k^2(w_k\cdot\nabla)w_k
  +J_k\bigl((B_k\cdot\nabla)w_k+(w_k\cdot\nabla)B_k\bigr).
\]
Subtract \eqref{eq:7.3} and divide by $J_k^2$ to obtain
\begin{equation}\label{eq:7.5}
 -\mu_k\Delta w_k+(w_k\cdot\nabla)w_k+\nabla p_k=e_k,
 \quad\operatorname{div}w_k=0,
 \quad\operatorname{Tr}w_k=J_k^{-1}a,
\end{equation}
where
\begin{equation}\label{eq:7.6}
 e_k=J_k^{-2}s_k
     -J_k^{-1}\bigl((B_k\cdot\nabla)w_k+(w_k\cdot\nabla)B_k\bigr).
\end{equation}
In particular $w_k-v_k=(A-B_k)/J_k\to0$ in $H^1$.
By $H^1\hookrightarrow L^6$ and \eqref{eq:7.4},
\begin{equation}\label{eq:7.7}
 \mu_k^{-1}\|e_k\|_{3/2}
 \le C_{\nu,\Omega}
 \bigl(J_k^{-1}\|s_k\|_\infty
        +\|w_k\|_{H^1}\|B_k\|_{H^1}\bigr)\to0.
\end{equation}
Consequently, uniformly over measurable subsets $D\subset\Omega$,
\begin{equation}\label{eq:7.8}
 \sup_{D\subset\Omega}
 \mu_k^{-1}\left|\int_D e_k\cdot w_k\right|
 \le\mu_k^{-1}\|e_k\|_{3/2}\|w_k\|_3=o(1).
\end{equation}
The exact cancellation of $(B_k\cdot\nabla)B_k$ in
\eqref{eq:7.6} will also be needed for an interior $L^2$ estimate.

\subsection{Pressure compactness and interior estimates}

Put $Y_k=\mu_k(w_k-A/J_k)$ and
$q_k=p_k-(\mu_k/J_k)\pi$. Then
\[
 Y_k\in V,\qquad Y_k\to0\text{ in }H^1,
 \qquad-\Delta Y_k+\nabla q_k
 =e_k-(w_k\cdot\nabla)w_k=:H_k.
\]
The hypotheses of \eqref{eq:3.5}--\eqref{eq:3.7} follow from
\[
 \sup_k\|H_k\|_{3/2}
 \le\sup_k\bigl(\|e_k\|_{3/2}
               +\|w_k\|_6\|\nabla w_k\|_2\bigr)<\infty
\]
and
\[
 \|H_k+\operatorname{div}(v\otimes v)\|_{H^{-1}}
 \le C\|e_k\|_{3/2}
   +(\|w_k\|_4+\|v\|_4)\|w_k-v\|_4\to0.
\]
Here $L^{3/2}\hookrightarrow H^{-1}$ is dual to
$H_0^1\hookrightarrow L^3$. On each fixed planar or off-axis
meridional Lipschitz patch, $H^1$ embeds into every finite
$L^q$. In the spatial case the spatial and meridional norms
are comparable away from the axis. Thus for $1<s<2$ and
$q=2s/(2-s)$,
\[
 \|\mu_k^{-1}e_k\|_{L^s}
 \le C\bigl(J_k^{-1}\|s_k\|_\infty
             +\|B_k\|_{H^1(\Omega)}\|w_k\|_{H^1(\Omega)}\bigr)\to0,
 \qquad\sup_k\|H_k\|_{L^s}<\infty
\]
on that patch. The constants may depend on the patch and on $s$.
The estimate \eqref{eq:3.9} gives
\[
 q_k\to p\text{ in }L^2,
 \qquad \nabla p=-(v\cdot\nabla)v.
\]
Hence $p_k\to p$ in $L^2$ and
$\Phi_k=p_k+|w_k|^2/2\to\Phi=p+|v|^2/2$ in $L^2$.
Lemma~\ref{lemma:5.1}, the calculation in \eqref{eq:5.9}
with lift $A$, and \eqref{eq:7.1} give
\[
 \operatorname{Tr}_{\Gamma_j}p=p_j,
 \qquad \nu=-\sum_jp_jF_j.
\]
Here the $F_j$ in the finite sum are the boundary fluxes
defined in \eqref{eq:1.flux}, not terms of the smooth force sequence.
Lemma~\ref{lemma:3.5} applied to $(Y_k,q_k)$ gives, on every
admissible component,
\[
 \lim_{\delta\downarrow0}\limsup_k
 \operatorname*{ess\,sup}_{0<h<\delta}
 \|q_k-p_j\|_{L^2(\Gamma_{j,h})}=0.
\]
The harmonic part $(\mu_k/J_k)\pi=(\mu_k^2/\nu)\pi$ satisfies
\eqref{eq:4.12}, while the kinetic estimate \eqref{eq:4.16}
uses $\operatorname{Tr}w_k=a/J_k$.
Thus \eqref{eq:4.19} and Lemma~\ref{lemma:4.3} apply to
$w_k,\Phi_k$, after the simultaneous shift \eqref{eq:5.10}
on zero-head components. The unscaled pressures $P_k$ retain
their mean-zero normalization in the final compactness argument.

\begin{lemma}[Interior control of the effective force]\label{lemma:7.1}
For every compact set $K\Subset\Omega$ in the plane, and every
compact rotational layer $K\Subset\Omega\cap\{r>0\}$ in space,
\[
 \|e_k\|_{L^2(K)}\longrightarrow0.
\]
\end{lemma}
\begin{proof}
Choose nested interior neighborhoods of $K$, separated from the
axis in space. The unit-viscosity equation for $\mu_kw_k$ has
right-hand side $H_k$ and pressure $p_k$. On the smaller
neighborhood $O$, interior Stokes regularity gives
\[
 \|\mu_kw_k\|_{W^{2,3/2}(O)}
 \le C_O\bigl(\|H_k\|_{3/2}
            +\|\mu_kw_k\|_{H^1}+\|p_k\|_2\bigr)\le C_O.
\]
Take rotational neighborhoods in space. On their generating
sections $r\in[r_*,r^*]$, the measures $2\pi r\,dr\,dz$ and
$dr\,dz$ are comparable. Cylindrical components have the same
meridional derivative bounds, while angular terms carry only
bounded factors $r^{-1}$ and $r^{-2}$.
In the plane, or on a smaller meridional patch, the
two-dimensional embeddings
$W^{2,3/2}\hookrightarrow W^{1,6}\hookrightarrow L^\infty$
therefore yield
\[
 \|\mu_kw_k\|_{L^\infty(K)}
 +\|\mu_k\nabla w_k\|_{L^6(K)}\le C_K.
\]
Using $1/2=1/3+1/6$ and $J_k^{-1}=\mu_k/\nu$ in
\eqref{eq:7.6}, we obtain
\begin{align*}
 \|e_k\|_{L^2(K)}
 &\le |\Omega|^{1/2}J_k^{-2}\|s_k\|_\infty\\
 &\quad+\nu^{-1}\bigl(
   \|B_k\|_3\|\mu_k\nabla w_k\|_{L^6(K)}
  +\|\nabla B_k\|_2\|\mu_kw_k\|_{L^\infty(K)}\bigr)
 \longrightarrow0
\end{align*}
by \eqref{eq:7.2} and \eqref{eq:7.4}.
\end{proof}

\subsection{The signed head identity for the corrected sequence}

Write $\omega_k=\operatorname{curl}w_k$, interpreted as scalar
vorticity in the plane. The smooth interior equations give
\begin{equation}\label{eq:7.9}
 \operatorname{div}(\nabla\Phi_k-e_k)
 =|\omega_k|^2+\mu_k^{-1}\operatorname{div}(\Phi_kw_k)
               -\mu_k^{-1}e_k\cdot w_k.
\end{equation}
Indeed, in space
$\nabla\Phi_k-e_k=w_k\times\omega_k
-\mu_k\operatorname{curl}\omega_k$. The planar formula follows
by adjoining a zero third component. The curl term has zero
total flux across each closed interior collar surface; Stokes'
formula uses only a neighborhood of that surface. Thus
\begin{equation}\label{eq:7.10}
 \left|\int_{\Gamma_h}(\nabla\Phi_k-e_k)\cdot n\,dS\right|
 \le C\int_{\Gamma_h}|w_k|\,|\nabla w_k|\,dS.
\end{equation}
This identity removes the need for a trace of the original
dual-space force on the collar.

\subsection{Boundary maximum and coarea for the force error}

Suppose \eqref{eq:5.24} holds and the maximal-pressure components
are admissible. Make the simultaneous shift in \eqref{eq:5.10}.
The Euler and local cycle hypotheses in
Lemmas~\ref{lemma:5.1} and \ref{lemma:5.2} follow from the
preceding estimates. The construction in
\eqref{eq:5.25}--\eqref{eq:5.29} supplies
$t_i\downarrow0$, $I_i=[5t_i/8,7t_i/8]$, and fixed interior
layers $L_i\Subset\Omega$ with $|L_i|\to0$; in space they are
separated from the axis. For $k\ge k(i)$ and almost every
$t\in I_i$, there are a free surface $S=S_{ik}(t)$ and a region
$D=D_{ik}(t,h)$ with
\[
 \partial D=S\mathbin{\dot\cup}\Gamma_h,
 \qquad \Phi_k|_S=-t,
 \qquad\nabla\Phi_k\cdot n_D=-|\nabla\Phi_k|\text{ on }S,
\]
\[
 S\subset L_i,\qquad \mathcal H^{n-1}(S)\ge c_*>0,
 \qquad \int_D|\omega_k|^2\ge c_i>0.
\]
Here $\Gamma_h$ is a collar of a fixed union $\Gamma$ of
zero-head physical components. The moving heads have the
local $W^{1,3/2}$ bound required for cycle convergence, by
the interior pressure estimate and
$\|w_k\|_6\|\nabla w_k\|_2\le C$.
The enstrophy lower bound uses weak convergence of $\omega_k$
and nonconstancy of the limiting head on a fixed interior
subregion. Use the Borel set $E_{ik}\subset L_i$ from
\eqref{eq:5.46}; its regular level sections are exactly $S_{ik}(t)$.

For fixed $i$, Lemma~\ref{lemma:4.3} and \eqref{eq:4.19}
give one common $h=h_k$, independent of $t$, such that
\[
 \int_{\Gamma_h}|w_k|\,|\nabla w_k|\le\epsilon,
 \qquad \|w_k\|_{L^2(\Gamma_h)}\le\sqrt{C_\epsilon}\,\mu_k,
 \qquad \|\Phi_k\|_{L^2(\Gamma_h)}\le\sigma.
\]
Incompressibility in $D$ and in the removed collar also gives
\[
 \mu_k^{-1}\int_S w_k\cdot n_D
 =-\nu^{-1}\sum_{\Gamma_j\subset\Gamma}F_j\le C_F,
 \qquad C_F=\nu^{-1}\sum_j|F_j|.
\]
All normals now point out of $D$. Integrating \eqref{eq:7.9}
and using the sign of the head gradient on $S$ yields
\begin{equation}\label{eq:7.11}
\begin{split}
 \int_S|\nabla\Phi_k|
 ={}&\int_{\Gamma_h}(\nabla\Phi_k-e_k)\cdot n_D
      -\int_S e_k\cdot n_D-\int_D|\omega_k|^2\\
 &+\frac t{\mu_k}\int_S w_k\cdot n_D
   -\frac1{\mu_k}\int_{\Gamma_h}\Phi_kw_k\cdot n_D
   +\frac1{\mu_k}\int_D e_k\cdot w_k.
\end{split}
\end{equation}
Choose first $\epsilon$ so that $C\epsilon<c_i/4$ in
\eqref{eq:7.10}; then choose $\sigma$ so that
$\sigma\sqrt{C_\epsilon}<c_i/4$. Next fix a sufficiently thin
collar and take $k$ sufficiently large for its bounds and for
the error in \eqref{eq:7.8} to be less than $c_i/4$.
This order of choices gives
\[
 \int_{S_{ik}(t)}|\nabla\Phi_k|
 \le C_Ft+\int_{S_{ik}(t)}|e_k|
 \quad\text{for almost every }t\in I_i.
\]
Put
\[
 X_{ik}=\int_{E_{ik}}|\nabla\Phi_k|^2,
 \qquad \eta_{ik}=\|e_k\|_{L^2(L_i)}.
\]
Coarea and Cauchy--Schwarz give
\[
 \int_{I_i}\int_{S_{ik}(t)}|e_k|\,dS\,dt
 =\int_{E_{ik}}|e_k|\,|\nabla\Phi_k|
 \le\eta_{ik}\sqrt{X_{ik}}.
\]
Consequently
\[
 X_{ik}\le\frac{3C_F}{16}t_i^2+\eta_{ik}\sqrt{X_{ik}},
 \qquad
 X_{ik}\le\frac{3C_F}{8}t_i^2+\eta_{ik}^2.
\]
For each fixed $i$, Lemma~\ref{lemma:7.1} permits choosing
$k=k(i)$ so large that $\eta_{ik}\le t_i$, in addition to all
collar and cycle requirements. Only a volume $L^2$ estimate
is used here, after integration in $t$; no trace of $f\in V'$
or of $F\in H^{-1}$ occurs. With $C$ independent of $i$,
\[
 \frac{c_*t_i}{4}
 \le\int_{E_{ik}}|\nabla\Phi_k|
 \le |L_i|^{1/2}\sqrt{X_{ik}}
 \le C|L_i|^{1/2}t_i.
\]
Division by $t_i$ contradicts $|L_i|\to0$.

\subsection{Interior and axial maxima}

If $\operatorname*{ess\,sup}_\Omega\Phi>\max_jp_j$, the
construction in \eqref{eq:5.51} gives
$D_{ik}(t)\Subset\Omega$ with its whole boundary equal to
$S_{ik}(t)$ and with enstrophy at least $c_i>0$.
Here $E_{ik}$ and $L_i$ denote the analogous interior level union
and separating layer, and $\eta_{ik}=\|e_k\|_{L^2(L_i)}$.
Their closures lie inside $\Omega$ and, in space, stay off the
axis, so Lemma~\ref{lemma:7.1} again applies for each fixed $i$.
Since $\Phi_k=-t$ on that boundary and $\operatorname{div}w_k=0$,
\[
 \int_{\partial D_{ik}(t)}\Phi_kw_k\cdot n_D=0.
\]
Equation~\eqref{eq:7.9}, estimate~\eqref{eq:7.8} and coarea
then give, for fixed $i$ and sufficiently large $k$,
\[
 X_{ik}+\frac12c_i|I_i|
 \le\eta_{ik}\sqrt{X_{ik}}
 \le X_{ik}+\frac14\eta_{ik}^2.
\]
Letting $k\to\infty$ contradicts $c_i|I_i|>0$.

It remains to consider a spatial boundary maximum attained on
an axis-intersecting component when the domain meets the axis.
All axis-intersecting components have the same Euler pressure.
Shift it to $p_0=0$, so $\Phi\le0$ almost everywhere. The other
$p_j\le0$, with at least one inner $p_j<0$ by
\eqref{eq:7.1} and $\sum_jF_j=0$. The limiting Euler equation
and $\operatorname{Tr}v=0$ give, as in
\eqref{eq:5.54}--\eqref{eq:5.57},
\[
 0\ge\int_\Omega(3\Phi-|v|^2/2)
 =\int_\Omega\operatorname{div}\bigl(xp+(x\cdot v)v\bigr)
 =-3\sum_{j=1}^Np_j|\Omega_j|>0.
\]
These alternatives exhaust the cases, since the essential
supremum of the head is at least every constant boundary trace.

\subsection{Passage to the limit and pressure representatives}

\begin{proof}[Completion of the proof of Theorem~\ref{thm:main}]
Every subsequence for which $J_k\to\infty$ has been excluded.
Hence $\sup_k\|U_k\|_{H^1}<\infty$ for the selected strongly
convergent sequence of representatives. Extract
\[
 U_k\rightharpoonup u\text{ in }H^1,
 \qquad U_k\to u\text{ in }L^4,
 \qquad\operatorname{Tr}u=a,
 \qquad\operatorname{div}u=0.
\]
Since
\[
 \|U_k\otimes U_k-u\otimes u\|_2
 \le(\|U_k\|_4+\|u\|_4)\|U_k-u\|_4\to0,
\]
passage to the limit gives
\[
 \nu(\nabla u,\nabla\eta)-(u\otimes u,\nabla\eta)
 =\langle F,\eta\rangle=\langle f,\eta\rangle
 \qquad(\eta\in V),
\]
which is \eqref{eq:1.weak}. Proposition~\ref{prop:quotient}
recovers $P_F\in L_0^2$ for this representative and, for the
same velocity, the pressure associated with every other
$H^{-1}$ extension of $f$. The fixed symmetry subspaces are
weakly closed, so all selected velocity symmetries pass to $u$.

In fact $d_k=U_k-u\in V$ is admissible in the difference of the
weak equations. With $T_k=U_k\otimes U_k-u\otimes u$,
\[
 \nu\|\nabla d_k\|_2^2
 =\langle F_k-F,d_k\rangle+(T_k,\nabla d_k),
\]
and therefore
\begin{equation}\label{eq:7.velocity-convergence}
 \nu\|\nabla d_k\|_2
 \le\|(F_k-F)|_V\|_{V'}+\|T_k\|_2\to0.
\end{equation}
For the chosen mean-zero pressures, the divergence inf--sup
estimate gives
\begin{equation}\label{eq:7.pressure-convergence}
 \|P_k-P_F\|_2
 \le C_\Omega\bigl(\nu\|\nabla d_k\|_2
                    +\|T_k\|_2+\|F_k-F\|_{H^{-1}}\bigr)\to0.
\end{equation}
Thus, along this subsequence,
\begin{equation}\label{eq:7.strong}
 (U_k,P_k)\longrightarrow(u,P_F)
 \quad\text{in }H^1(\Omega;\mathbb R^n)\times L_0^2(\Omega).
\end{equation}
The pressure estimate uses convergence in the full $H^{-1}$
space, which is available because the representative $F$ was
fixed before approximation.
\end{proof}

\begin{remark}[What solenoidal convergence does and does not control]
\label{rem:pressureconvergence}
The velocity estimate \eqref{eq:7.velocity-convergence} only
uses the $V'$ norm of the force difference. The analogous
claim for unadjusted pressures is false. For a nonconstant
$\psi\in C_c^\infty(\Omega)$, take
\[
 a=0,\qquad U_k=0,\qquad G_k=k\nabla\psi.
\]
Then $G_k|_V=0$ for all $k$, while the associated mean-zero
pressures are $k(\psi-\psi_\Omega)$ and diverge in $L^2$.
Thus a mere substitution of $V'$ for $H^{-1}$ in
\eqref{eq:7.pressure-convergence} would be incorrect.
If $g_k=G_k|_V\to f$ in $V'$ and specific representatives
$G_k\in H^{-1}$ are used, write
\[
 G_k-\mathcal Eg_k=\nabla q_k,\qquad q_k\in L_0^2.
\]
Whenever the corresponding velocities converge as in
\eqref{eq:7.velocity-convergence}, the adjusted pressures
$P_{G_k}-q_k$ converge in $L^2$ to $P_{\mathcal Ef}$, because
$\mathcal Eg_k\to\mathcal Ef$ in $H^{-1}$.
\end{remark}

\begin{remark}[The natural energy-dual endpoint]\label{rem:endpoint}
For every $H^1$ velocity satisfying the solenoidal weak equation,
\[
 |\langle f,\eta\rangle|
 \le\bigl(\nu\|\nabla u\|_2+\|u\|_4^2\bigr)\|\nabla\eta\|_2
 \qquad(\eta\in V).
\]
Thus the effective force of a finite-energy solution necessarily
belongs to $V'$. This explains the force space in the main
statement and the absence of any additional regularity
requirement on an independently chosen gradient part.
\end{remark}

\subsection*{Conflict of interest statement}
The author states that there is no conflict of interest.

\subsection*{Data availability statement}
Data sharing is not applicable to this article as no datasets were generated or analysed for this study.

\subsection*{AI declaration}
The author acknowledges the use of generative AI tools in developing and checking mathematical arguments and in drafting and revising the
manuscript. The author takes full responsibility for the mathematical content.



\begin{thebibliography}{99}

\bibitem{BKK}
J.~Bourgain, M.~Korobkov and J.~Kristensen,
On the Morse--Sard property and level sets of Sobolev and BV functions,
\emph{Rev. Mat. Iberoam.} \textbf{29} (2013), 1--23.
\href{https://doi.org/10.4171/RMI/710}{doi:10.4171/RMI/710}.

\bibitem{Galdi}
G.~P.~Galdi,
\emph{An Introduction to the Mathematical Theory of the Navier--Stokes
Equations: Steady-State Problems}, second edition, Springer, 2011.
\href{https://doi.org/10.1007/978-0-387-09620-9}{doi:10.1007/978-0-387-09620-9}.

\bibitem{GQ}
J.~Gopalakrishnan and W.~Qiu,
Partial expansion of a Lipschitz domain and some applications,
\emph{Front. Math. China} \textbf{7} (2012), 249--272.
\href{https://doi.org/10.1007/s11464-012-0189-2}{doi:10.1007/s11464-012-0189-2}.

\bibitem{HMT}
S.~Hofmann, M.~Mitrea and M.~Taylor,
Geometric and transformational properties of Lipschitz domains,
Semmes--Kenig--Toro domains, and other classes of finite perimeter domains,
\emph{J. Geom. Anal.} \textbf{17} (2007), 593--647.
\href{https://doi.org/10.1007/BF02937431}{doi:10.1007/BF02937431}.

\bibitem{KLW2022}
M.~Korobkov, W.~Lyu and S.~Weng,
On the existence of helical invariant solutions to steady
Navier--Stokes equations,
\emph{SIAM J. Math. Anal.} \textbf{54} (2022), no.~3, 3974--3996.
\href{https://doi.org/10.1137/20M1343324}{doi:10.1137/20M1343324}.

\bibitem{KPR2013}
M.~V.~Korobkov, K.~Pileckas and R.~Russo,
On the flux problem in the theory of steady Navier--Stokes equations
with nonhomogeneous boundary conditions,
\emph{Arch. Ration. Mech. Anal.} \textbf{207} (2013), 185--213.
\href{https://doi.org/10.1007/s00205-012-0563-y}{doi:10.1007/s00205-012-0563-y}.

\bibitem{KPR2015axis}
M.~V.~Korobkov, K.~Pileckas and R.~Russo,
An existence theorem for steady Navier--Stokes equations in the
axially symmetric case,
\emph{Ann. Sc. Norm. Super. Pisa Cl. Sci.} (5) \textbf{14} (2015), 233--262.
\href{https://doi.org/10.2422/2036-2145.201204_003}{doi:10.2422/2036-2145.201204\_003}.

\bibitem{KPR2015}
M.~V.~Korobkov, K.~Pileckas and R.~Russo,
Solution of Leray's problem for stationary Navier--Stokes equations
in plane and axially symmetric spatial domains,
\emph{Ann. of Math.} (2) \textbf{181} (2015), 769--807.
\href{https://doi.org/10.4007/annals.2015.181.2.7}{doi:10.4007/annals.2015.181.2.7}.

\bibitem{KPRbook}
M.~Korobkov, K.~Pileckas and R.~Russo,
\emph{The Steady Navier--Stokes System:
Basics of the Theory and the Leray Problem},
Advances in Mathematical Fluid Mechanics,
Birkh\"auser, Cham, 2024.
\href{https://doi.org/10.1007/978-3-031-50898-1}{doi:10.1007/978-3-031-50898-1}.

\bibitem{KY2009}
H.~Kozono and T.~Yanagisawa,
Nonhomogeneous boundary value problems for stationary Navier--Stokes
equations in a multiply connected bounded domain,
\emph{Pacific J. Math.} \textbf{243} (2009), 127--150.
\href{https://doi.org/10.2140/pjm.2009.243.127}{doi:10.2140/pjm.2009.243.127}.

\bibitem{MW2012}
M.~Mitrea and M.~Wright,
\emph{Boundary Value Problems for the Stokes System in Arbitrary
Lipschitz Domains}, Ast\'erisque \textbf{344} (2012), viii+241 pp.
\href{https://doi.org/10.24033/ast.808}{doi:10.24033/ast.808}.
The regularity theorem cited here is Theorem~1.1.3 of the published volume.

\end{thebibliography}
\end{document}